\documentstyle{amlts}
\begin{document}
\annalsline{155}{2002}
\received{October 16, 2000}
\startingpage{235}
\def\bye{\end{document}}
 \font\tenrm=cmr10

%--------------- Author macros ---------------
%for Bbb in amstex
\catcode`\@=11
\font\twelvemsb=msbm10 scaled 1100
\font\tenmsb=msbm10
%\font\ninemsb=msbm7 scaled 1100%msbm9
\font\ninemsb=msbm10 scaled 800
\newfam\msbfam
\textfont\msbfam=\twelvemsb  \scriptfont\msbfam=\ninemsb
  \scriptscriptfont\msbfam=\ninemsb
\def\msb@{\hexnumber@\msbfam}
\def\Bbb{\relax\ifmmode\let\next\Bbb@\else
 \def\next{\errmessage{Use \string\Bbb\space only in math
mode}}\fi\next}
\def\Bbb@#1{{\Bbb@@{#1}}}
\def\Bbb@@#1{\fam\msbfam#1}
\catcode`\@=12

 \catcode`\@=11
\font\twelveeuf=eufm10 scaled 1100
\font\teneuf=eufm10
\font\nineeuf=eufm7 scaled 1100%eufm9
\newfam\euffam
\textfont\euffam=\twelveeuf  \scriptfont\euffam=\teneuf
  \scriptscriptfont\euffam=\nineeuf
\def\euf@{\hexnumber@\euffam}
\def\frak{\relax\ifmmode\let\next\frak@\else
 \def\next{\errmessage{Use \string\frak\space only in math
mode}}\fi\next}
\def\frak@#1{{\frak@@{#1}}}
\def\frak@@#1{\fam\euffam#1}
\catcode`\@=12
%-------------- Author entries --------------------

%\intro %(Optional, Introduction)

\newcommand{\be}{\begin{equation}}
\newcommand{\ee}{\end{equation}}
\newcommand{\ba}{\begin{array}}
\newcommand{\ea}{\end{array}}
\newcommand{\bea}{\begin{eqnarray}}
\newcommand{\eea}{\end{eqnarray}}
\newcommand{\bee}{\begin{eqnarray*}}
\newcommand{\eee}{\end{eqnarray*}}

%le petit carre de fin de demonstration
\def \trait (#1) (#2) (#3){\vrule width #1pt height #2pt depth #3pt}
\def \fin{\null\hfill
        \trait (0.1) (5) (0)
        \trait (5) (0.1) (0)
        \kern-5pt
        \trait (5) (5) (-4.9)
        \trait (0.1) (5) (0)
\medskip}

\def\R{{\bf R}}
\def\N{{\bf N}}

\def\lim{\mathop{\rm lim}}
\def\goto{\rightarrow}
\def\supp{{\rm supp}~}
\def\sup{\mathop{\rm sup}}
\def\exp{{\rm exp}}
\def\e{\varepsilon}
\def\w{w}
\def\l{\lambda}
\def\ch{{\rm ch}}
\def\sh{{\rm sh}}
\def\th{{\rm th}}
\def\ln{{\rm ln}}
\def\ath{{\rm ath}}
\def\atan{{\rm atan}}
\def\ds{\displaystyle}
\def\ts{\textstyle}
\def\om{\omega}
\def\Ai{{\rm Ai}}
\def\vp{\varphi}
\def\s{\sigma}
\def\q{H}
\def\qb{\overline q}
\def\wb{\overline \w}
\def\a{\alpha}
\def\ab{\overline \alpha}
\def\b{\beta}
\def\bb{\overline \beta}
\def\g{\gamma}
\def\gb{\overline \gamma}
\def\db{\overline \delta}
\def\d{\delta}
\def\W{W}
\def\H{H}
\def\Lun{L_1}
\def\Hb{\overline H}
\def\Lunb{\overline L_1}
\def\QXX{{Q\over 2}+yQ_y}
\def\ind{{\rm ind}}
\def\zz{|}
\def\etun{\eta_{I,n}}
\def\etde{\eta_{II,n}}
\def\etu{\eta_{I}}
\def\wun{w_{1,n}}
\def\wdeux{w_{2,n}}
\def\deg{{\left({Q\over 2}+yQ_y\right)}}

%%%%%%%%%%%%%%%%%%%%
%%%%%%%%%%%%%%%%%%%%
%%%%%%%%%%%%%%%%%%%%
%%%%%%%%%%%%%%%%%%%%
%%%%%%%%%%%%%%%%%%%%

 \title{Stability of  blow-up profile\\ and lower bounds for blow-up rate\\
for the critical generalized KdV  equation}
\shorttitle{The critical generalized KdV equation} 
  \twoauthors{Yvan Martel}{Frank Merle}
 \institutions{Universit\'e de Cergy-Pontoise, Cergy-Pontoise, France\\
{\eightpoint {\it E-mail address\/}: Yvan.Martel@math.u-cergy.fr}\\
\vglue6pt
Universit\'e de Cergy-Pontoise, Cergy-Pontoise, France\\
Institut Universitaire de France, Paris,  France\\
{\eightpoint {\it E-mail address\/}: Frank.Merle@math.u-cergy.fr}}
%-------------- Article Text--------------------

\centerline{\bf Abstract}
\vglue12pt
The generalized Korteweg-de Vries equations are a class of Hamiltonian
systems in infinite dimension derived from the KdV equation where the\break 
quadratic term is replaced by a higher order power term. These equations
have two conservation laws in the energy space $H^1$ ($L^2$ norm and
energy). We consider  in this paper the {\it critical} generalized 
KdV equation, which corresponds to the smallest power of the nonlinearity 
such that the two conservation laws do not imply 
a bound in $H^1$ uniform in time for all $H^1$ solutions
(and thus global existence).  

From \cite{M3}, there do exist for this equation solutions $u(t)$
such that\break $|u(t)|_{H^1}\goto +\infty$ as $t\uparrow T$, where $T\le +\infty$
(we call them blow-up solutions). The question is to describe, in a 
qualitative way,  how blow up occurs.

For solutions with $L^2$ mass close to the minimal mass allowing blow up
and with  
decay in $L^2$ at the right, we prove after rescaling and translation
which leave invariant the $L^2$ norm that the solution converges 
to a {\it universal} profile locally in space at the blow-up time
$T$.
From the nature of this profile, we improve the standard lower bound
on the blow-up rate for finite time blow-up solutions.

\section{Introduction}

1.1.
We consider in this paper the critical generalized Korteweg-de Vries equation:
\be
\left\{
\ba{ll}
u_t+(u_{xx}+u^5)_x=0 , &\quad (t,x)\in \R^+\times \R,\\[5pt]
u(0,x)= u_0(x) ,&\quad x\in \R,
\ea
\right.\label{kdv}
\ee
for $u_0 \in H^1(\R)$.
It is a special case of   the generalized Korteweg-de Vries equations, $p\ge 2$ integer:
\be
\left\{
\ba{ll}
u_t+(u_{xx}+u^p)_x=0 , &\quad (t,x)\in \R^+\times \R,\\[5pt]
u(0,x)= u_0(x) ,&\quad x\in \R.
\ea
\right.\label{kdvp}
\ee
The case $p=2$ corresponds to the KdV equation (see Korteweg and de Vries \cite{KDV}), 
and $p=3$ to the modified KdV equation.
These  two cases have been studied extensively for being completely
integrable (see, for example, Lax \cite{LAX}
and Miura \cite{MIURA}).
For all $p$, these equations are considered as universal models of Hamiltonian systems.
From this Hamiltonian structure, there are two conservation laws
\begin{eqnarray}
\label{masseintro}\int u^2(t)&\hskip-7pt=\hskip-7pt&\int u_0^2   \hbox{ (mass conservation),}\\[4pt]
 \label{energieintro} 
\qquad {1\over 2} \int u^2_x(t)-{1\over p+1} \int u^{p+1}(t)&\hskip-7pt=\hskip-7pt&
{1\over 2} \int u^2_{0x} -  {1\over p+1} \int u^{p+1}_0\\&&\hskip.95in
 \hbox{ (energy conservation).}\nonumber\end{eqnarray}

From \cite{KPV}, we have the following existence and uniqueness
result in the energy space $H^1(\R)$: 
for $u_0\in H^1(\R)$, there exists $T>0$ and
 a unique maximal solution 
$u\in C([0,T),H^1(\R))$ 
of (\ref{kdvp}) on $[0,T)$. Moreover, either $T=+\infty$, 
or $T<+\infty$, and then $|u(t)|_{H^1}\goto +\infty$, as
$t\uparrow T$. In addition,
$\hbox{for all } t\in [0,T)$, (\ref{masseintro}) and (\ref{energieintro}) are satisfied.
For equation (\ref{kdv}), the local Cauchy problem is also well
posed in $L^2(\R)$ (see \cite{KPV}). We refer to Kato \cite{K} and Ginibre and Tsutsumi \cite{GT} for 
previous results on the well-posedness of the Cauchy problem for (\ref{kdvp}). See Bourgain \cite{B}
for the periodic case.

For $p<5$ (the subcritical case), as a consequence of the Gagliardo-Nirenberg inequality,
all solutions in $H^1$ are global and bounded in time.

In this paper, we consider only the critical case $p=5$. We define the energy
$$E(u)={1\over 2} \int u^2_x -{1\over 6} \int u^6.$$
We consider  solutions in the energy space $H^1(\R)$. 
Let us introduce the ground state $Q$, unique positive solution (up to translation) of 
$$Q_{xx}+Q^5=Q,\quad Q\in H^1(\R),\quad  \hbox{or equivalently}\quad
Q(x)={3^{1/4}\over \ch^{1/2}(2x)}.$$
Note that $u(t,x)=Q(x-t)$ is a special solution of (\ref{kdv}), and $E(Q)=0$.

On the one hand, the variational  characterization of $Q$; i.e.: for $v\in H^1(\R)$,
\begin{eqnarray}\label{charac}
&&\hbox{if  } 0<\int v^2\le \int Q^2 \hbox{ and } E(v)\le 0,
\\
&&\hbox{then there exists } \l_0>0, x_0\in \R ~ / ~ v=\l_0^{1/2} Q(\l_0 (x-x_0)),
\nonumber\end{eqnarray}
which provides the following
Gagliardo-Nirenberg inequality with \pagebreak best constant (see Weinstein \cite{W2}):
\be\label{gn}
  \hbox{for all } v\in H^1(\R),\quad
   {1\over 6}\int v^6 \le {1\over 2} 
          \left( {\int v^2 \over \int Q^2} \right)^2 \int v_x^2,  
\ee
implies that for $|u_0|_{L^2}<|Q|_{L^2}$, the solution $u(t)$ is global and uniformly bounded in $H^1$.

On the other hand, for $|u_0|_{L^2}>|Q|_{L^2}$ there is no obstruction to blow up from 
energy-type arguments.
Existence of solutions  of (\ref{kdv}) blowing up in finite or infinite time in the energy space  $H^1$ 
has been  proved by Merle (see \cite{M3} and also Martel and Merle \cite{MM2}). 
More precisely: \vglue-24pt
\phantom{HI}
\begin{quote} There exists $\alpha_0>0$ such that for all $u_0\in H^1(\R)$, if
$E(u_0)<0$ and $\int u_0^2 <\int  Q^2 +\alpha_0$,   then the solution $u(t)$ blows up in $H^1$ in finite or
infinite time.
\end{quote}
\vglue-9pt\noindent
Note that numerical observations  suggest  existence of blow up in
finite time; see Bona et al.\ \cite{Bona}.
The argument for the blow-up proof in \cite{M3} 
is not direct. Arguing by contradiction, we consider a limit object, recurrent
in time, as  $t\goto +\infty$. The idea is to show that the recurrence in time yields some rigidity on this
object. Then, we are able to prove both  elliptic and oscillatory integral type estimates on this limit
solution; together with the three conservation laws (mass, energy and the additional invariant
$\int u(t,x) dx$ when $u_0\in L^1$), this gives a contradiction with the Liouville Theorem in [12].

Moreover, the quantity $|Q|_{L^2}$ is the minimal amount of $L^2$ norm 
that concentrates at blow-up time, in the sense that
for some function $x(t)$, we have 
 for all $\e_0>0$,  $\underline{\lim}_{t\uparrow T} \int_{|x-x(t)|\le \e_0} |u(t)|^2 \ge \int Q^2.$ 
\vglue3pt
In the study of the blow-up phenomenon, there are two main questions:    the profile at blow-up time 
(in  some smaller scale, describing the blow-up dynamics)
and the rate of blow up of the solutions. These two questions are clearly linked.
We consider the question of the blow-up profile for initial data as in the blow-up result, i.e.\
with $\int u_0^2\le \int Q^2 +\alpha_0.$

The first result of this paper is a characterization of the blow-up profile, which is $Q$, up to the invariances
of the equation. This is, in some sense, a generalization of the Liouville theorem   in \cite{MM2}, and of its corollary
which says that any bounded solution, starting close to $Q$ in $H^1$, converges locally in space to $Q$
for large time.

\proclaimtitle{Stability of $Q$ as a blow-up profile}
\specialnumber{1}\proclaim{Theorem}\label{THPROFILE}
There exists $\alpha_0>0$ such that if $u_0\in H^1(\R)$ satisfies
$$\int u_0^2 <\int Q^2 +\alpha_0$$
and if the solution $u(t)$ of {\rm (\ref{kdv})} blows up in finite or infinite   time $T>0${\rm ,} then
for all $0\le t<T${\rm ,} there exists $\l(t)>0$ and $x(t)\in \R$ such that either\pagebreak
\begin{eqnarray*}
\l^{1/2}(t) u(t,\l(t)x+x(t))\rightharpoonup Q && \hbox{as $t\uparrow T$ in
$H^1(\R)$ weak},\\ \noalign{\noindent or} 
 -\l^{1/2}(t) u(t,\l(t)x+x(t))\rightharpoonup Q && \hbox{as $t\uparrow T$ in
$H^1(\R)$ weak.}\\ \noalign{\vskip-28pt}
\end{eqnarray*}
\endproclaim 

{\it Remark}. Note that the alternative in Theorem \ref{THPROFILE} 
comes from the fact that if $u(t,x)$ is a solution of (\ref{kdv}) then 
$-u(t,x)$ is also a solution of (\ref{kdv}). 

There is no such result of determination of blow-up
profile in Hamiltonian systems or more generally in evolution partial differential equations except
for diffusion equations where the existence of 
Liapunov functions plays a fundamental role. Indeed, in the case of the nonlinear heat equations
$u_t=\Delta u +u^p$ in $\R^N$, with some restriction on $p$, the blow-up rate and 
profile (and their stability)
have been determined; see for example Giga and Kohn \cite{GK} and Fermanian-Kammerer, Merle and
Zaag \cite{FMZ}.

\demo{{R}emark} If $E(u_0)<0$ in Theorem \ref{THPROFILE}, we actually know that $u(t)$ blows up in finite
or infinite time; see \cite{M3}. Note that the result implies the stability and the universality of the blow-up 
profile in the region $\int u_0^2\le \int Q^2 +\alpha_0$ and $E(u_0)<0$.
\enddemo  

We consider now an application of this result to obtain a refined lower bound on the blow-up rate.
In particular we exclude some candidates (deduced from a  scaling argument) of the blow-up rate.
By a scaling argument and the resolution of the Cauchy problem, 
if $u(t)$ is a solution blowing up at some finite time $T>0$, then for some $C>0$,
$$\hbox{for all } t_0\in [0,T),  \quad  |u_x(t_0)|_{L^2}\ge {C\over (T-t_0)^{1/3}}.$$
Indeed, consider
$$v_{t_0}(t,x)=|u_x(t_0)|_{L^2}^{-1/2} u(t_0+|u_x(t_0)|_{L^2}^{-3}t,|u_x(t_0)|_{L^2}^{-1}x);$$
$v_{t_0}$ is a solution of (\ref{kdv}) by scaling invariance.
We have $|v_{t_0x}|_{L^2}+|v_{t_0}|_{L^2}\le C$, and so by the resolution of the Cauchy problem
 locally in time by a fixed-point argument
(see \cite{KPV}),
there exists $\tau>0$, independent of $t_0$, such that $v_{t_0}(t)$ is defined on $[0,\tau]$.
Therefore, $t_0+|u_x(t_0)|_{L^2}^{-3} \tau<T$, which is the desired result. Theorem 1 implies that  this
lower bound represents the exact blow-up rate for no solution with small $L^2$ mass. 
Indeed, we have the following theorem.
 
\proclaimtitle{Lower bound on the blow-up rate}
\specialnumber{2}\proclaim{Theorem}\label{THRATE}
There exists $\alpha_0>0$ such that if $u_0\in H^1(\R)$ satisfies
$$\int u_0^2 <\int Q^2 +\alpha_0$$
and if the solution  $u(t)$ of {\rm (\ref{kdv})} blows up in finite time $T>0${\rm ,} then
$$\lim_{t\uparrow T} (T-t)^{1/3}|u_x(t)|_{L^2}=+\infty.$$
\endproclaim

\vglue-20pt
{\it Remark}.
For the critical KdV equation, Bona and Weissler \cite{BW} constructed explicit solutions
of (\ref{kdv}) with self similar blow up. However these solutions are not in the physical space $L^2$ and also
exist in the linear context. In fact,
Theorem 2 excludes this type of blow-up rate  for solutions with $L^2$ mass close to the minimal mass allowing blow
up,  and we expect this result to extend to all initial data in 
$H^1$. 

This phenomenon, which forces the solutions to blow up more quickly than the self similar rate in the
energy space, seems to be typical of Hamiltonian systems in PDE with infinite speed of propagation.
It is still an open problem for the critical nonlinear Schr\"odinger equation (NLSE):
\be
\left\{
\ba{l}
{\rm i} \,u_t=-u_{xx}-|u|^4u , \quad (t,x)\in \R^+\times \R,\\[5pt]
u(0,x)= u_0(x) ,\quad x\in \R.
\ea
\right.\label{nlse}
\ee
Indeed, the rate predicted by scaling arguments in energy space, which is $|u_x(t)|_{L^2}
\sim C/\sqrt{T-t}$, should be not relevant for blow-up solution in $H^1$. 
Evidence  for this fact is on the one hand the existence of explicit
solutions of the nonlinear critical Schr\"odinger equation with a  blow-up rate  $C/(T-t)$, and
on the other hand numerical results suggesting different blow-up rates. See for example Merle \cite{M2}
for  more information about blow up for NLSE.
Note that for the Zakharov system (critical NLSE coupled with a wave equation), 
the optimal lower bound  for the blow-up rate ($C/(T-t)$) has been derived by Merle \cite{MerleZakharov}. 

\vglue2pt

The second author thanks Stanford University  where part of this work was done.

\demo{{\rm 1.2.} Strategy of the proof}
First, as in \cite{MM2} and \cite{M3}, from the fact that $\int u_0^2<\int Q^2+\hskip1pt\alpha_0$,
with $\alpha_0$ small, a parametrization of the problem allows us to see the evolution in time
of the size and location of the solution. Note that  for all $\l_0>0, x_0\in \R$, 
$v_{\l_0,x_0}(x)=\l_0^{1/2} Q(\l_0(x -x_0))$ is such that $\int v_{\l_0,x_0}^2=\int Q^2$ and
$E(v_{\l_0,x_0})=0$. Therefore, there is no obstruction from the conservation laws to the existence of  a solution of the form
$$u(t,x)\sim \l^{-1/2}(t) Q(\l^{-1}(t) (x -x(t))).$$
In fact, it follows from variational arguments that a blow-up solution with $\alpha_0$ small is close
in $H^1$ to the set $\{v_{\l_0,x_0}\}$, for all time  close to the blow-up time. Indeed, set  
\be\label{defe} 
  \e(t,y)=\lambda^{1/2}(t)u(t,\lambda(t)y+x(t))
    -Q(y),
\ee for $\lambda(t)>0$, $x(t)$, two $C^1$ functions
to be chosen later. For a suitable choice of $\l(t)$ and $x(t)$, $\e(t)$ is 
uniformly small in time in $H^1(\R)$. Change 
 the time variable as follows:
\be
  s=\int_0^t {dt' \over \l^3(t')},
\quad \hbox{or equivalently,}\quad {ds\over dt}={1\over
\l^3};
\ee
then $\e(s)$ satisfies, for $s\ge 0$, $y\in \R$, 
\begin{eqnarray}
&& \label{eqn400}\\
 \ds  \e_s \ds& =& (L\e)_y   + {\l_s \over \l} \left({Q\over 2} 
      +yQ_y\right) + \left( {x_s \over \l}-1 \right) Q_y\nonumber \\[4pt]
       &&  +\ {\l_s \over \l} \left({\e\over 2} 
      +y\e_y\right) + \left( {x_s \over \l}-1 \right) \e_y
         -(10Q^3\e^2+10Q^2\e^3+5Q \e^4+\e^5)_y,\nonumber
\end{eqnarray}
where 
\be
L \e =-\e_{xx}+\e-5 Q^4 \e 
    =-\e_{xx}+\e-{15 \over \ch^2(2x)}\e.
\ee
(See Lemma 1 in \cite{MM1}.)
Note that $x(t)$ and $\l(t)$ are geometrical parameters
related to the two invariances of equation (\ref{kdv}):
respectively, translation and dilation invariances.

If, for all $t\ge 0$, 
$u(t)$  is sufficiently close to $Q$ in $H^1$, up to scaling and translation, 
we can  define a unique $C^1$ function $s\goto (\l(s),x(s))$ such that for all
$s\ge 0$, 
$$
\int y\left({Q\over 2}+yQ_y\right) \e(s)=\int y Q_y \e(s)=0.$$
The reason to consider such  orthogonality conditions 
on $\e(s)$ is the fact that
they are adapted  to a  Virial-type
identity on $\e(s)$ (${d\over ds}\int y\e^2(s)$; see \S 2.2). Indeed, these relations
cancel some interactions in the Virial relation, and are one of the crucial tools of this paper.

  Theorem 2 follows directly from Theorem 1, and the
fact that $u(t,x)=Q(x-t)$ is a solution of (\ref{kdv}) such that for all $t$, $|u_x(t)|_{L^2}=|Q_x|_{L^2}$.
(See \S 4.1.) The proof of Theorem 1 is by contradiction. Assume that we have a solution for $\alpha_0$ small such that
$$\e(s)\not \rightharpoonup 0  \quad \hbox{in $H^1(\R)$, as $s\goto +\infty$}.$$
As in \cite{M3}, the idea is to define a recurrent object as $t\goto T$. From the property
of recurrence of this object, and some  almost monotonic in-time functional, this object has more
properties (decay properties as $y\goto -\infty$).
The proof of nonexistence of such an object then concludes the proof of
Theorem 1.

More precisely, define $s_n\goto +\infty$ such that $\e(s_n)\rightharpoonup \widetilde \e(0)\not \equiv 0$,
and $\widetilde u$ solution of (\ref{kdv}) with $\widetilde u(0)=Q+\widetilde \e(0)$. This  limit
solution
$\widetilde u(t)$  is associated to $\widetilde \e$, $\widetilde \l$,
$\widetilde x$, and is such that $\widetilde \l(0)=1$. Define $\tau>0$ such that for all $ s\in [0,\tau)$, $
{1\over 1.1} \le \widetilde \l(t)\le 1$.
The contradiction follows from three facts.
\smallbreak
(i) We first prove exponential decay on the left for $\widetilde \e$, in the sense that 
$$\hbox{for all  } s\in [0,\tau) ,   \hbox{ for all }y<0 , \qquad
\quad 
\widetilde \e(s,y)|\le C(\alpha_0) e^{-{|y|\over 12}},$$
where $C(\alpha_0)\goto 0$ as $\alpha_0\goto 0$. (See \S 2.3 and \S 4.2.)

We then conclude by using two dispersion relations giving information on the dynamics of $\widetilde \e(s)$:
\smallbreak
(ii) an $L^1$ relation, involving a quantity of the type $\int \widetilde \e(s) V$, where $V$ is bounded such that 
$V(y)\goto 1$ as $y\goto -\infty$, and $V(y)\goto 0$ as $y \goto +\infty$.
\smallbreak
(iii) an $L^2$ relation, which is a local Virial-type identity, i.e.\ an expression for
${d\over ds}\int \Psi_A \widetilde \e^2,$ where $\Psi_A$ is bounded and such that  $\Psi_A(y)\sim y$ for $|y|<A$.
\smallbreak
These three facts provide a contradiction on $\widetilde \e$ and $\widetilde u$, in the two possible
regimes $\tau<+\infty$ or $\tau=+\infty$. (See \S 3.2 and \S 4.2.)

\vglue4pt

The paper is organized as follows. In Section 2, we establish fundamental relations. 
In Section~3, we study a simpler and
more geometric case, where the limit object is more natural that in Theorem 1.
This allows us to present the main ideas without the technical difficulties of Section~4. We prove directly that for
$\alpha_0$ small, the upper bound
$$|u_x(t)|_{L^2}\le {C\over (T-t)^{1/3}}$$
is not possible for blow-up solutions. Section 4 is then devoted to the proof of Theorems 1 and 2.

\vglue-8pt
\section{Energy and dispersive relations}
\advance\eqcount by 11
\vglue-4pt

2.1. {\it Decomposition of the solution}.
For $v\in H^1(\R)$, let
$$\alpha(v)=\int v^2-\int Q^2,\quad E(v)={1\over 2}\int v_x^2 -{1\over 6}\int v^6.$$

For $\alpha_0>0$ to be fixed later, let $u_0\in H^1(\R)$ be such that
$\alpha(u_0)\le \alpha_0$. Assume that the solution $u(t,x)$ of (\ref{kdv}) blows up
in finite or infinite time $T>0$ (which implies by (\ref{gn}) that $\alpha(u_0)\ge 0$).
Note that  $-u(t,x)$ is also a solution of~(\ref{kdv}) with initial data $-u_0$, and
with the same properties as $u(t,x)$.

First, we have the following decomposition of the solution, using variational tools, and 
elementary geometrical properties. This decomposition allows us to extract some finite dimensional
approximation of the solution which gives in some sense the time evolution of its size and its location.

\proclaimtitle{Sharp decomposition and modulation of the solution}
\specialnumber{1} \proclaim{Lemma} \label{MODULATION} 
  There exists $\alpha_1>0$ such that if
$\alpha_0<\alpha_1$ then there exist $0\le t(u_0)< T$ and  continuous functions $\l:[t(u_0),T)\goto (0,+\infty),$
$x:[t(u_0),T)\goto \R${\rm ,} such that{\rm ,} for $v\equiv u$ or $v\equiv -u${\rm ,} for all $t\in [t(u_0),T)${\rm ,}
\be\label{defofeps}    \e(t,y)=\l^{1/2}(t) v(t,\l(t)y+x(t))-Q(y)
\ee
satisfies the following properties\/{\rm :} for all $t\in [t(u_0),T),$
\be\label{propofeps1}  
\int yQ_y \e(t)=\int y\left({\ts {{Q\over 2}+yQ_y}}\right) \e(t)=0,
\ee
\be\label{propofeps2}
\left|1-\l(t) {|u_x(t)|_{L^2}\over |Q_x|_{L^2} }\right|
+|\e(t)|_{H^1}\le \delta'(\alpha_0),\quad \hbox{where $\delta'(\alpha)\goto 0$ as $\alpha\goto 0$.}
\ee
\endproclaim

\vglue-20pt
{\it Remark}. In the rest of this paper, we assume that $v=u$ in Lemma \ref{MODULATION}, by possibly replacing
$u(t,x)$ by $-u(t,x)$.

\medskip

\demo{{P}roof} The proof is similar to the one of Lemmas 1 and 2 in \cite{M3}.
\enddemo

{\it Step} 1. 
We first  prove that
any blow-up solution with $L^2$ norm close to the $L^2$ norm of $Q$ is close to the set
$\{\pm \l_0^{-1/2} Q(\l_0^{-1}(x-x_0))\mid\l_0>0, x_0\in \R \}$ for $t$ close enough to the blow-up time.

\nonumproclaim{Claim}  There exists  $\alpha_2>0$ such that
the  following property 
holds true{\rm .} For all $0<\alpha'\le \alpha_2${\rm ,} there exists $\delta=\delta(\alpha')> 0${\rm ,} with
$\delta(\alpha')\goto 0$ as $\alpha'\goto 0${\rm ,} such that  for all $u\in H^1(\R)${\rm ,} $u\not \equiv 0${\rm ,} if
$$\alpha(u)\le \alpha',\quad {E(u)\le \alpha' \int u_x^2},$$
then there exist  $ x_0\in \R$ and $\epsilon_0\in \{-1,1\}$ such that
$$|Q-\epsilon_0 \l_0^{1/2}u(\l_0 x+x_0 )|_{H^1}\le \delta(\alpha'),$$
with $\l_0={|Q_x|_{L^2}/|u_x|_{L^2}}.$
\endproclaim

\demo{{P}roof of the claim}
It follows from variational arguments.
For the sake of contradiction, consider a sequence $(u_n)$ of functions in $H^1(\R)$, $u_n\not
\equiv 0$, such that 
$$\lim_{n\goto \infty} \int u^2_n\le \int Q^2 ,\quad \lim_{n\goto \infty} {E(u_n)\over \int u_{nx}^2}\le 0.$$

Let us recall the following variational result, following the variational characterization of $Q$
(see Lemma 1 in \cite{M3}):
if a sequence $(v_n)$ of $H^1$ functions satisfies:
\be\label{444}
\lim_{n\goto +\infty} \int v_n^2= \int Q^2,\quad \int v_{nx}^2=\int Q_x^2,\quad \mathop{\overline{\lim}}_{n\goto
+\infty} E(v_n)\le 
0,\ee
then there exists a sequence $(x_n)$ of $\R$, and $\epsilon_0\in \{-1,1\}$ such that
\be\label{claimlemma1}\ds \lim_{n\goto \infty} \epsilon_0 v_n(.+x_n)=Q \hbox{ in $H^1(\R)$ as $n\goto \infty$}.
\ee

We set
$$\l_n={|Q_x|_{L^2} \over |u_{nx}|_{L^2}} \quad \hbox{and}\quad 
v_n=\l_n^{1/2} u_n(\l_n x).$$  
Note that $\int v_n^2=\int u_n^2$ and $\int (v_n)^2_x=\int Q_x^2.$
We prove that the sequence $v_n$ satisfies (\ref{444}), which finishes the proof of the claim.
\enddemo
 \pagebreak

Indeed, by Gagliardo-Nirenberg inequality with best constant (\ref{gn}), we have
$${E(u_n)\over \int u_{nx}^2 }  \ge {1\over 2} \left(1-{\int u_n^2 \over \int Q^2}\right),$$
and since $\lim_{n\goto +\infty} \alpha(u_n)\le 0$, $\lim_{n\goto +\infty} {E(u_n)\over \int u_{nx}^2}\le 0$,
 it follows that
$$\lim_{n\goto \infty} \int u_n^2 =\lim_{n\goto \infty} \int v_n^2=\int Q^2\quad \hbox{and}
\quad \lim_{n\goto \infty} {E(u_n)\over \int u_{nx}^2}=0.$$
By direct calculations, we have:
$$E(v_n)={E(u_n)\over \int u_{nx}^2} \int Q_x^2,\quad 
\int (v_n)_x^2=\int Q_x^2.$$
Therefore, the sequence $(v_n)$ satisfies (\ref{444}).

\demo{Step 2. {\it Modulation of the solution}}
\quad 
By conservation of energy, we have  for all $t\in [0,T)$, $E(u(t))=E(u_0)$. Since $\lim_{t\uparrow T} |u_x(t)|_{L^2}
=+\infty$, there exists $t_0\in [0,T)$ such that 
$$\hbox{for all } t\in [t_0,T[, \quad  
\alpha_0 \int u_x^2(t)\ge E(u(t))=E(u_0).$$
 Therefore, by step 1, for all 
$t\in [t_0,T)$, there exist $x_0(t)\in \R$, and $\epsilon_0(t)\in
\{-1,1\}$ such that, with $\l_0(t)=|Q_x|_{L^2}/|u_x(t)|_{L^2},$
$$|Q-\epsilon_0(t) \l_0^{1/2}(t) u(t,\l_0(t)x+x_0(t))|_{H^1}\le \delta(\alpha_0).$$
As in Merle \cite{M3}, if $\delta(\alpha_0)<|Q|_{L^2}/4$ (which is true when
$\alpha_0$ is  small), then $\epsilon_0(t)$ is  independent of $t$.
Moreover, we can assume, with no restriction, that $\epsilon_0(t)=1$
(note that $-u(t,x)$ is also a solution of (1)).

Now, as in Lemma 2 of \cite{M3}, we sharpen the decomposition; i.e.\ we choose $\l(t)$, $x(t)$ close to $\l_0(t)$, $x_0(t)$
such that $\e(t)=\l^{1/2}(t) u(t,\l(t)y+x(t))-Q(y)$  is small in $H^1$ and also satisfies
suitable  orthogonality conditions:
 for all $t\in [t(u_0),T)$, 
$$
\int yQ_y \e(t)=\int y\left({\ts {{Q\over 2}+yQ_y}}\right) \e(t)=0.
$$
Here, we use the implicit function theorem.
 Note that with respect to Lemma 2 in \cite{M3}, we have modified the orthogonality
conditions. See also Part B of \cite{MM2}.
In the present case, 
we have
 \begin{eqnarray*}&&\hskip-36pt
  \left({d\over d \l_0} \l_0^{1/2} Q(\l_0 x+x_0)\right)_{(\l_0=1,x_0=0)}
 \\[3pt]
 &&\qquad
 ={Q\over 2}+y Q_y, 
\left({d\over d x_0} \l_0^{1/2} Q(\l_0 x+x_0)\right)_{(\l_0=1,x_0=0)}=Q_y,  
 \end{eqnarray*}
\pagebreak

\noindent
and the nondegeneracy conditions are satisfied since
$$\int \left({\ts {Q\over 2}}+y Q_y\right) y \left({\ts {Q\over 2}}+y Q_y\right)=0,\quad 
\int (Q_y) y \left({\ts {Q\over 2}}+y Q_y\right)=\int \left({\ts {Q\over 2}}+y Q_y\right)^2\not = 0,$$
$$\int \left({\ts {Q\over 2}}+y Q_y\right) yQ_y =\int \left({\ts {Q\over 2}}+y Q_y\right)^2\not =0,\quad 
\int (Q_y) yQ_y=0,$$
where we have used $\int Q \left({\ts {Q\over 2}}+y Q_y\right)=0$, and parity properties.
\hfill\qed
\enddemo

Note that by invariance of the equation by translation in time, we may assume that $t(u_0)=0$ in the rest
of this paper. We also assume that $u(t,x)$ satisfies the decomposition of Lemma \ref{MODULATION}, and
not $-u(t,x)$. This is the first possibility in Theorem 1.

Let 
\be\label{time} s=\int_0^t {dt'\over \l^3(t')}, \quad \hbox{or equivalently}, \quad
{ds\over dt}={1\over \l^3}.\ee

Observe that by the scaling property of equation (\ref{kdv}) and local well-posedness
of (\ref{kdv}) in $H^1$ (see introduction), we have  by (\ref{propofeps2})
\be\label{lower}
\hbox{for all } t\in [0,T),\quad \l(t)\le C (T-t)^{1/3}.
\ee
Therefore, when $t$ takes its values on  $[0,T)$, $s$ takes its values in all $\R^+$.
In the next lemma, following \cite{MM1}, and \cite{M3}, Lemma 3 and Corollary 1, we gather useful properties of $\e(s)$,
$\l(s)$ and $x(s)$.

\proclaimtitle{Properties of the
decomposition}
\specialnumber{2} \proclaim{Lemma} \label{PROPERTIESDECOMPOSITION} There exists $0<\alpha_3\le \alpha_2${\rm ,}
 such that if
$\alpha_0<\alpha_3$ then  $\l(s)$ and $x(s)$ are $C^1$ functions on $\R^+$ and  the following
properties exist\/{\rm :} 
\vglue6pt
{\rm (i)} Equation of $\e(s)${\rm .}  The function $\e(s)$ satisfies{\rm ,} for $s\in \R^+, y\in \R,$
\bee  \e_s & = &(L\e)_y   + {\l_s \over \l} \left({Q\over 2} 
      +yQ_y\right) + \left( {x_s \over \l}-1 \right) Q_y\\
&&  +\ {\l_s \over \l} \left({\e\over 2} 
      +y\e_y\right) + \left( {x_s \over \l}-1 \right) \e_y 
         -(R(\e))_y, 
\eee
where 
$
L \e =-\e_{xx}+\e-5 Q^4 \e$ and $R(\e)=10Q^3\e^2+10 Q^2\e^3+5 Q\e^4+\e^5.$\vglue6pt
{\rm (ii)} Smallness properties{\rm .} For some $C>0${\rm ,}
\be\label{smallness}\hbox{for all } s\ge 0,\quad |\e(s)|_{L^2}+|\e_y(s)|_{L^2}\le C\sqrt{\alpha_0},
\ee
\be\label{17bis}\hbox{for all } t\in [0,T),\quad { |u_x(t)|_{L^2}\over 2|Q_x|_{L^2}}
\le {1\over \l(t)}\le {2 |u_x(t)|_{L^2}\over |Q_x|_{L^2}}.
\ee
\vglue6pt
{\rm (iii)} Control of the geometrical parameters{\rm .} For some $C>0${\rm ,} 
\be\label{corollaire1} \hbox{for all } s\ge 0,\quad 
\left| {\l_s\over \l}\right| +\left|{x_s\over \l}-1\right|  \le C \left(\int \e^2 e^{-{|y|\over 2}}\right)^{1/2}
\le C \sqrt{\alpha_0}.
\ee
\endproclaim 

\vglue-20pt
{\it Proof}.
For (i) see  \cite{MM1}, and \cite{M3}.
Let us   
recall the following structural properties of $L$:
\be L(Q^3)=-8 Q^3,\quad L(Q_y)=0,\ee
\be\label{coerd}
\hbox{for all } \e\in H^1(\R),\hbox{ if } \int Q^3 \e=\int Q_y \e=0 \hbox{ then } (L\e,\e)\ge \int \e^2.\qquad 
\ee

 We prove (ii) using some
ideas from \cite{M3}.
By the definition of $\e$ and the conservation of mass, we have
$$2\int Q\e+\int \e^2=\int u^2-\int Q^2=\int u_0^2-\int Q^2=\alpha_0.$$
By the conservation of energy, we have
$\l^2 E(Q+\e)=E(u(t))=E_0<0$, and by direct calculations, 
\bea \label{preener}
&& E(Q+\e)+\left(\int Q\e+{1\over 2} \int \e^2 \right)
\\
&&\qquad =\ {1\over 2} (L\e,\e)  -{1\over 6} \left[ 20\int Q^3\e^3 +15 \int Q^2 \e^4 +6\int Q \e^5+\int \e^6
\right]. \nonumber
\eea
Therefore, 
\be\label{firstener}
(L\e,\e)\le \alpha_0 + C |\e|_{H^1} |\e|_{L^2}^2.
\ee
Note that by the choice of  orthogonality conditions on $\e$,
this is not sufficient to conclude the proof of (ii) directly. Indeed, they are suitable for the Virial identity but not for
the energy identity. Nevertheless, 
 consider an auxiliary function
$\e_1= \e -a \deg -b Q_y,$ where
$$\int \e_1 Q^3=\int \e_1 Q_y=0,$$
 and take
$a\!=\!{\int \e Q^3 \over \int \deg Q^3},$ $
b\!=\!{\int \e Q_y \over \int Q_y^2}$
(note that $\int \deg Q^3 \!=\!{1\over 4}\int Q^4\! \not \equiv\!0$). 
Note that we  also have
$\e=\e_1+a \deg + bQ_y,$ so that by orthogonality conditions on $\e$,
$
a={\int \e_1 \deg\over \int \deg^2},$
$b={\int \e_1 y \deg \over \int \deg^2}.$
\vglue4pt
Now, since $\int \deg Q=0$, $L\deg=-2 Q$, $LQ_y=0$, we find after some
elementary calculations:
$$\int \e Q=\int \e_1 Q,\quad (L\e,\e)=(L\e_1,\e_1)- 4 a (\e,Q).$$
By the expressions for $a$ and $b$, we have for some constant $K$,
\be\label{norml2}
{1\over K} (\e_1,\e_1)\le (\e,\e)\le K (\e_1,\e_1).
\ee
Thus, from (\ref{coerd}) and (\ref{firstener}), 
$${1\over K} (\e,\e)\le (\e_1,\e_1)\le (L\e_1,\e_1) \le \alpha_0+4|a| |(\e,Q)|
+C |\e|_{L^2}^2 |\e|_{H^1}.$$
For $\alpha_0$ small, $|\e|_{H^1}$ and $|a|$ are  small and from the conservation of mass, we 
have $2 |(\e,Q)|\le \alpha_0+\int \e^2$; thus 
$${1\over K} (\e,\e)\le 2 \alpha_0 + C|\e|_{L^2}^2 (|a|+|\e|_{H^1})
\le 2\alpha_0+{1\over 2K}|\e|_{L^2}^2.$$
Therefore, $(\e,\e)\le 4 K  \alpha_0$ and from (\ref{firstener}),
$(L\e,\e)\le C\alpha_0.$
The conclusion then comes from the fact that
$|\e|^2_{H^1}\le (L\e,\e) + 5 \int Q^4 \e^2 
\le (L\e,\e)+c (\e,\e).$

For (iii)
note that by multiplying the equation of $\e(s)$ by $yQ_y$ 
and then by $y\left({\ts {Q\over 2}}+yQ_y\right)$, 
using the decay property of $Q$ at infinity, we obtain
\bee
&&\hskip-36pt\left| {\l_s\over \l} \left(\mu_0 +\int ({\ts {\e\over 2}}+y\e_y) yQ_y)\right)
-\int \e L((yQ_y)_y) \right|\\
&&\qquad \le\ C \left|{x_s\over \l}-1\right| \left(\int \e^2 e^{-{|y|\over 2}}\right)^{1/2} +C \int
\e^2e^{-{|y|\over 2}},
\eee
and
\bee 
&&\hskip-36pt\left| \left({x_s\over \l}-1\right) \left(\mu_0 +\int  \e_y y({\ts{Q\over 2}}+ yQ_y )\right)
-\int \e L((y({\ts{Q\over 2}}+ yQ_y ))_y) \right|\\
& &\qquad\le\  C \left| {\l_s\over \l}
\right| \left(\int \e^2 e^{-{|y|\over 2}}\right)^{1/2}  + C \int \e^2e^{-{|y|\over 2}},
\eee
where $\mu_0=\int ({\ts{Q\over 2}}+ yQ_y )^2>0.$
Therefore, for $\alpha_3$ small enough,
\bee &&\hskip-36pt
\left| {\l_s\over \l} -  {1\over \mu_0} \int \e L((yQ_y)_y) \right|
+\left| \left({x_s\over \l}-1\right) - {1\over \mu_0} \int \e L((y({\ts{Q\over 2}}+yQ_y))_y)\right| \\
&&\quad \le C \left(\left| {\l_s\over \l}
\right|+\left|{x_s\over \l}-1\right|
\right) \left(\int \e^2 e^{-{|y|\over 2}}\right)^{1/2} + C \int \e^2e^{-{|y|\over 2}},
\eee
and (\ref{corollaire1}) follows.
Note that we have in addition
\begin{eqnarray*}
\left| {\l_s\over \l} -  {1\over \mu_0} \int \e L((yQ_y)_y) \right|
&+& \left| \left({x_s\over \l}-1\right) - {1\over \mu_0} \int \e L((y({\ts{Q\over 2}}+yQ_y))_y)\right| \\[4pt]
& \le&
C\int \e^2(s) e^{-{|y|\over 2}}. \\
\noalign{\vskip-36pt}
\end{eqnarray*}
\hfill\qed
\vglue12pt

In the rest of this paper, $u(t)$ denotes a solution blowing up in finite or infinite time $T>0$ such that 
$\int u^2= \int Q^2+\alpha_0$, with $\alpha_0<\alpha_3$; we consider the decomposition in $\e$, $\l$ and $x$,
satisfying the properties given in Lemmas \ref{MODULATION} and~\ref{PROPERTIESDECOMPOSITION}.

\demo{{\rm 2.2.} Relations}
We consider four different quantities for $\e(s)$. Two of them are related to conservation laws, and the others
concern dispersion relations respectively in $L^1$ and in $L^2$.

\proclaimtitle{Mass and energy conservation}
\specialnumber{3} \proclaim{Lemma} \label{MASSENERGY}
\be\label{masse}
2\int \e Q +\int \e^2=\alpha_0 \quad \hbox{{\rm (}\/mass conservation\/{\rm ),}}
\ee
\be\label{energie}
\left| \l^2 E_0 + \int \e Q -{\ts {1\over 2}} \int \e_y^2 \right|
\le C \int \e^2 e^{-|y|} +C \alpha_0^2 \int \e^2_y \quad \hbox{{\rm (}\/energy relation\/{\rm )}}.
\ee
\endproclaim

\demo{Proof} 
 First, we recall that the conservation of the $L^2$ norm of $u(t)$  and the notation
$\alpha_0=\int u^2-\int Q^2$ give (\ref{masse}) directly.

Second, by (\ref{preener}), 
the decay properties of $Q$, $|\e(s)|_{L^\infty}\le C\sqrt{\alpha_0}$ 
(by (\ref{smallness})),
and the Gagliardo-Nirenberg
inequality $\int \e^6 \le C(\int \e^2)^2 \int \e_y^2 \le C\alpha_0^2 
\int \e_y^2$,  we obtain (\ref{energie}). \enddemo

Next,   defining 
$$J(s)=\left[\int \e(s,y) \left(\int_y^{+\infty} \left({Q\over 2} +z Q_z\right) dz \right)dy\right]- {1\over 4}
 \left(\int Q\right)^2,$$
we have the following lemma.

\proclaimtitle{$L^1$-type dispersion}
\specialnumber{1}\proclaim{Proposition}\label{FORMULEJ}
Assume that $J(s_0)$ is well-defined for some $s_0$. Then $J(s)$ is well\/{\rm -}\/defined for all $s\ge 0$ 
and is of class $C^1${\rm .} Moreover{\rm ,}
\be\label{formuleJ}
\left| J_s+{\l_s\over 2\l} J +2\int \e Q \right|\le C\int \e^2 e^{-{|y|\over 2}}.
\ee
\endproclaim

\demo{Proof} We obtain the equation satisfied by $J_s$ and $J$ by multiplying the equation of $\e$
(see Lemma \ref{PROPERTIESDECOMPOSITION} (i))
by $\int_{y}^{+\infty} {Q\over 2}+zQ_z$ and then integrating by parts. This calculation is formal but can be justified
rigorously by regularization arguments; see the proof of Lemma 6 in \cite{MM1}. We obtain
\bee
J_s+{\l_s\over 2\l}J+2\int \e Q&=&{\l_s\over \l} \int y\left({Q\over 2}+yQ_y\right) \e
+\left({x_s\over \l}-1\right) \int \left({Q\over 2}+yQ_y\right) \e \\ &&  -\ \int 
\left({Q\over 2}+yQ_y\right)R(\e).\eee
We have used in particular
$L\left({Q\over 2}+yQ_y\right)=-2Q,$ $\int Q\left({Q\over 2}+yQ_y\right)=0$ and
$$\int \left({Q\over 2}+yQ_y\right) \int_y^{+\infty} \left({Q\over 2}+zQ_z\right)={1\over 2}\left(\int{Q\over 2}+yQ_y\right)^2
={1\over 8} \left(\int Q\right)^2.$$

From Lemma \ref{PROPERTIESDECOMPOSITION}, we have  
$$\left|\int R(\e)
({\ts {Q\over 2}}+yQ_y)\right| + \left|{\l_s\over \l}\right|^2
+\left|  {x_s\over \l}-1\right|^2\le C \int \e^2 e^{-{|y|\over 2}},$$
and
$$\left(\int \left({Q\over 2}+yQ_y\right) \e\right)^2\le C \left( \int (|y| e^{-(1-{1\over \sqrt{2}}) |y|})
(|\e|e^{-{|y|\over \sqrt{2}}})\right)^{2}\le C \int \e^2 e^{-{|y|\over 2}},
$$
which proves (\ref{formuleJ}).\enddemo

Finally, we will need a Virial type relation. In \cite{MM1}, we   considered a quantity of the type
$\int y\e^2$. Unfortunately, in our situation, we do not have   good control on $\e$ on the right
(i.e.\ as $y\goto +\infty$).
Therefore, we need a localized version of this identity, which is given in the next
lemma.

Let $\Phi\in {\cal C}^2(\R)$, $\Phi(x)=\Phi(-x)$, $\Phi'\le 0$ on $\R^+$, 
such that 
$$\Phi(x)=1 \hbox{ on $[0,1]$;}\quad \Phi(x)=e^{-x} \hbox{ on $[2,+\infty)$,}\quad
e^{-x}\le \Phi(x)\le 3 e^{-x} \enspace \hbox{on
$\R^+$.}$$
Let $$\Psi(x)=\int_0^x \Phi(y) dy.$$ Note that $\Psi$ is an odd function,
$\Psi(x)=x$ on $[-1,1]$, and $|\Psi(x)|\le 3$ on $\R$.

For a parameter $A>0$, we set
$$\Psi_A(x)=A \Psi\left({x\over A}\right),\quad
\hbox{so that}\quad \Psi'_A(x)=\Phi\left( {x\over A}\right)=\Phi_A(x), \quad \hbox{and}$$
$$\Psi_A(x)=x \hbox{ on $[-A,A]$},\quad |\Psi_A(x)|\le 3A \hbox{ on $\R$},
\quad e^{-{|x|\over A}} \le \Phi_A(x)\le 3 e^{-{|x|\over A}} \hbox{ on $\R$.}
$$

\proclaimtitle{Local Virial relation, $L^2$-type dispersion}
\specialnumber{2}\proclaim{Proposition}\label{VirialLOCAL}
 xThere exists $A_0>2${\rm ,} $\alpha_5>0$ and  $\delta_0>0$ such that for $\alpha_0<\alpha_5${\rm ,}
\begin{eqnarray}
&&\hbox{if for all } s\ge 0,\quad \int yQ_y \e(s) =\int y\left({Q\over 2}+yQ_y\right)\e(s) =0, \enspace
\hbox{then}\label{Viriallocal}\\
&& \left(\int \Psi_{A_0} \e^2\right)_s
 \le -\delta_0 \int (\e^2+\e_y^2)e^{-{|y|\over A_0}}
+{1\over \delta_0} \left(\int \e Q \right)^2.\nonumber
\end{eqnarray}
\endproclaim

\demo{Proof} First we note that this Virial relation is reminiscent of
the Virial identity for $\e$  established in \cite{MM1}.
Indeed, when $\int y \e^2$ is defined, we have from \cite[Lemma 5]{MM1}:
\bee
&& \hskip-48pt \left({1\over 2} \int y \e^2\right)_s+{\l_s\over \l} {1\over 2} \int y \e^2 \\
&&\quad =\
     {\l_s \over \l} \int y\left({Q\over 2}+yQ_y\right)\e
       +  \left( {x_s\over \l}-1 \right) \left( \int yQ_y\e
	-{1\over 2}\int \e^2\right) \\
       & & \qquad    -\ {3\over 2} (L\e,\e) +\int \e^2 
              -10 \int Q^3 \left({Q\over 2}+yQ_y\right) \e^2 \\
       & & \qquad  
           +\ 10\int \left({2Q^3\over 3}-yQ^2Q_y\right) \e^3
	+5\int \left({3Q^2\over 2}-yQ Q_y\right) \e^4 \\
& &  \qquad 	+\ \int(4Q-yQ_y)\e^5+{5\over 6}\int \e^6.
\eee
This Virial identity is obtained formally by multiplying the equation of $\e$ by
$y \e$, integrating  by parts, and  using some properties of $Q$ and the operator $L$.
See Lemma 5 in \cite{MM1} for a complete proof and regularizing argument.
The formula above does not require any orthogonality conditions on $\e(s)$.

In \cite{MM2}, we   established the following.
Denote
\bee
H_{\infty}(\e,\e)&=& {3\over 2} (L\e,\e)-\int \e^2 +10 \int Q^3 \left({Q\over 2}+y Q_y\right)\e^2 \\
&=& {3\over 2} \int \e_y^2 +{1\over 2}\int \e^2-{5\over 2} \int Q^4 \e^2 +10\int y Q_y Q^3 \e^2.
\eee
There exists $\delta_1>0$ such that 
$$\hbox{if}\quad \int Q \e=\int y\left({Q\over 2}+yQ_y\right) \e=0 \quad
\hbox{then} \quad H_{\infty}(\e,\e)\ge \delta_1 \int \left(\e_y^2 +\e^2\right).$$
(See Proposition 4 in \cite{MM2}, and its proof. Note that this result is optimal.)

Therefore, if 
$$\int yQ_y \e=\int Q \e=\int y\left({Q\over 2}+yQ_y\right) \e=0,$$
then by the Virial identity and the control on ${\l_s\over \l}$, ${x_s\over \l}-1$ and $|\e|_{L^\infty}$
given by Lemma \ref{PROPERTIESDECOMPOSITION} (ii), (iii), we obtain
$$
\left( \int y \e^2\right)_s+{\l_s\over \l}  \int y \e^2
\le -2 \delta_1 \int \left(\e_y^2 +\e^2\right) +C\left(\int \e^2\right)^{3/2}
\le -\delta_1 \int \left(\e_y^2 +\e^2\right) ,$$
for $|\e|_{H^1}$ small.

In fact,
the orthogonality condition $\int Q\e=0$ can be replaced by an additional term in the 
Virial relation :
if only $\int yQ_y \e=\int y\left({Q\over 2}+yQ_y\right) \e=0$ then
$$\left( \int y \e^2\right)_s+{\l_s\over \l}  \int y \e^2
\le  -{\delta_1\over 2} \int \left(\e_y^2 +\e^2\right) + C \left(\int Q \e\right)^2,$$
for $|\e|_{H^1}$ small.

The proof of Proposition~2 is a local version of this result, seen as a perturbation of the previous
identity 
(note that $\Psi_A$ is a cut-off of $y$ since $\Psi_A(y)=y $ on $[-A,A]$). 
Since the arguments are rather technical, we present the rest of the proof in
Appendix A. \enddemo
\pagebreak

\demo{{\rm 2.3.} Exponential decay on the left}
Together with relations given in the previous subsections, we now introduce a fundamental tool
which links $L^2$ concentration as a Dirac mass at the blow-up time 
(or $L^2$ compactness for large time) with exponential decay on the left.
This technique was introduced in \cite{M3}.
\enddemo

\proclaimtitle{Exponential decay in $L^2$ norm on the left}
\specialnumber{4} \proclaim{Lemma} \label{firstexpo}
There exists $\alpha_6>0${\rm ,} $C_1${\rm ,} $C_2>0$ satisfying the following property{\rm .}
Suppose $\alpha_0<\alpha_6${\rm .}
Suppose that $u^2(t)$ blows up as a Dirac mass at the blow\/{\rm -}\/up time{\rm ,} in the sense that
$$u^2(t,x+x(t))\rightharpoonup\left(\int u_0^2\right) \delta_{x=0},\quad \hbox{as $t\uparrow T$}.$$
If $s_1<s_2$ satisfy
$$\hbox{for all } s\ge s_1,\enspace \l(s)\le \l(s_1),\enspace \hbox{and for all } s\in [s_1,s_2],\enspace {\l(s_1)\over 1.1}\le
\l(s)\le
\l(s_1),$$ then
$$\hbox{for all } y_0<0,\hbox{ for all } s\in [s_1,s_2],\quad 
\int_{y<y_0} \e^2(s) \le  C_1 e^{-C_2 |y_0|}.$$
\endproclaim

Note that as in \cite{M3}, we need a control on the oscillations of $|u_x(t)|_{L^2}$\break (or $\l(t)$) which allows blow up.
As a direct consequence, we have:

\proclaimtitle{Exponential decay in $L^\infty$ norm on the left}
\specialnumber{1} \proclaim{{C}orollary}\label{secondexpo}
Under the  assumptions of the preceding lemma{\rm ,}  
$$\hbox{for all } y<0,\hbox{ for all } s\in [s_1,s_2],\quad
|\e(s,y)|\le C'_1 \alpha_0^{1/4} e^{-C'_2 |y|}.$$
\endproclaim

{\it Proof of Corollary {\rm \ref{secondexpo}}}. 
From the proof of Gagliardo-Nirenberg inequality, and Lemmas \ref{PROPERTIESDECOMPOSITION} and
\ref{firstexpo}, we have:   
  for all $y_0<0$  and  $s\in [s_1,s_2]$, 
\vglue9pt
\hfill ${\displaystyle |\e(s,y_0)|\le \left(\int \e_y^2\right)^{1/4} \left(\int_{y<y_0} \e^2\right)^{1/4}
\le C'_1 \alpha_0^{1/4} e^{-C'_2 |y_0|}.}$ \hfill\qed\vglue12pt

 The proof of Lemma \ref{firstexpo} follows essentially from techniques
introduced in \cite{M3} and \cite{MM2}.
We give some notation, and then we recall a result from \cite{M3} concerning the solutions
of (\ref{kdv}).

Let $K=2\sqrt{3}$ (for example). Define
$$\hbox{for all } x\in\R,\enspace \phi(x)=cQ\left({x\over K}\right),\enspace
\psi(x)=\int_{-\infty}^x \phi(y) dy,\enspace {\rm where}\enspace c={1\over K\int_{-\infty}^{+\infty} Q},$$
so that
$$\hbox{for all } x\in\R,\quad 0\le \psi(x)\le 1, \lim_{x\goto -\infty} \psi(x)=0,\quad
\lim_{x\goto +\infty} \psi(x)=1.$$

Consider $z(t)$ a solution of (\ref{kdv}) satisfying the decomposition of \S 2.1, with parameters
$\l_z$ and $x_z$.
For $x_0\in\R$, and $t_0\ge 0$ define
$$\hbox{for all } t\ge 0,\quad
{\cal I}_{x_0,t_0}(t)=\int z^2(t,x) \psi\left(x-x_z(t_0)-x_0-{\ts {1\over 4}}(x_z(t)-x_z(t_0))\right) dx.
$$

\proclaimtitle{Almost monotonicity of the mass on the left
\cite{M3}}
\specialnumber{5} \proclaim{Lemma} \label{monotonicite} Suppose that  for all $t\ge t_0${\rm ,}
 $0<\l_z(t)\le 1.1${\rm .} There exist $\alpha_6>0${\rm ,}
$C_0>0$  and $a_0>0$ such that if $\alpha_0<\alpha_6${\rm ,} then
$$\hbox{for all } x_0\le -a_0,\hbox{ for all } t\ge t_0,\enspace
{\cal I}_{x_0,t_0}(t)-{\cal I}_{x_0,t_0}(t_0)\le C_0 e^{x_0 \over 3}.$$
\endproclaim

{\it Sketch of the proof of Lemma {\rm \ref{monotonicite}}}.
Let $\widetilde x=x-x_z(t_0)-x_0-{\ts {1\over 4}}(x_z(t)-x_z(t_0))$.
We have by direct calculations
\begin{eqnarray*}
{\cal I}'_{x_0,t_0}(t)&=&-3 \int z_x^2(t,x) \phi(\widetilde x)-{(x_z)_t \over 4} \int z^2(t,x)
\phi(\widetilde x) \\
&&+\ \int z^2(t,x)\phi''(\widetilde x)+{5\over 3} \int z^6(t,x)\phi(\widetilde x),
\end{eqnarray*}
and
$\phi''\le {1\over K^2} \phi,$ ${9\over 10 \l_z^2} \le
(x_z)_t.$
It follows from 
  control of the last term of the expression of ${\cal I}'_{x_0,t_0}(t)$, using the   decomposition of the solution $z$, that 
\be\label{plusprecis}
 {\cal I}'_{x_0,t_0}(t)\le C (x_z)_t e^{-{3\over 4K} (x_z(t)-x_z(t_0))+{x_0\over K}}.
\ee
Then the result follows by integration in time. See \cite{M3} for details. \hfill\qed 

\demo{Proof of Lemma {\rm \ref{firstexpo}}} 
Let $s_1<s_2$ be as in the statement of Lemma \ref{firstexpo}, and let $s_0\in [s_1,s_2]$.
We apply Lemma \ref{monotonicite} to a rescaled version of $u(t)$.

Let $\l_0=\l(s_0)$, and
$$z(t',x')=\l_0^{1/2} u(\l_0^3 t',\l_0 x').$$
The function $z$ is a solution  of (\ref{kdv}).
We note that $\l_z(s)=\l(s)/\l_0$ and $x_z(s)=x(s)/\l_0$  are the geometrical parameters associated to $z$.
Now
$$z^2(t',x'+x_z(t'))\rightharpoonup \left(\int u^2_0 \right)\delta_{x'=0},\quad \hbox{as $t'\uparrow  T/\l_0^3$,}$$
$$\hbox{for all } s\ge s_0,\quad \l_z(s)={\l(s)\over \l(s_0)}\le {\l(s_1)\over \l(s_0)}\le 1.1,$$
since $s_0\in [s_1,s_2]$.
Let $y_0<-a_0$ and let $s\ge s_0$.
We apply Lemma \ref{monotonicite} to $z$ between $s_0$ and $s\ge s_0$. We obtain
\be\label{aujour}\begin{array}{ll} & \ds
\int z^2(s_0,x')\psi(x'-x_z(s_0)-y_0) dx'\\ &\ds
\qquad 
\ge -C e^{y_0\over 3} +\int z^2(s,x')\psi\left(x'-x_z(s_0)-y_0-{\ts{1\over 4}}(x_z(s)-x_z(s_0))\right)dx'.
\end{array}
\ee

Let  $v(s,y)=\e(s,y)+Q(y)=\l^{1/2}(s) u(s,\l(s) y + x(s))$,
so that $v(s_0,y)=z(s_0,y+x_z(s_0));$
the left-hand side of (\ref{aujour}) is
$$\int v^2(s_0,y) \psi(y-y_0) dy.$$

On the other hand, since $x_z(s)\ge x_z(s_0)$
(by (\ref{corollaire1}), ${d x_z\over ds}\ge 0$), 
and $\psi$ is nondecreasing, we have
$$\int v^2(s_0,y) \psi(y-y_0) dy\ge  -C e^{y_0\over 3}+ \int z^2(s,x')\psi(x'-x_z(s)-y_0)dx', $$
which gives, as $s\goto +\infty$ (which corresponds to $t'\goto T/\l_0^3$) ,
$$\int v^2(s_0,y) \psi(y-y_0) dy \ge -C e^{y_0\over 3} +
\left(\int u_0^2 \right) \psi(-y_0).$$

Now we use the properties of the function $\psi$ and $\int v^2(s_0)=\int u_0^2$.
Therefore,
\bee
 \int_{y<y_0} v^2(s_0,y) dy & \le &2 \int_{y<y_0} v^2(s_0,y) \left(1-\psi(y-y_0)\right) dy \\
&\le& 2 \int v^2(s_0,y) \left(1-\psi(y-y_0)\right) dy\\
&  \le& 2\left[ \int u_0^2 + Ce^{y_0\over 3} -\psi(-y_0) \int u_0^2 \right]\\
&\le& C  [e^{y_0\over 3}+ (1-\psi(-y_0))] \le C e^{c' y_0}.
\eee
The proof of the lemma then follows from $\e(s_1,y)=v(s_1,y)-Q(y)$, and the
decay properties of $Q$. 
Note that $|\e(s)|_{L^\infty}\le C \sqrt{\alpha_0}$ (by Gagliardo-Nirenberg inequality
and Lemma \ref{PROPERTIESDECOMPOSITION})
extends the inequality for all $y<0$.
\enddemo

\section{Direct proof of a lower bound at blow-up time}
\advance\eqcount by 32

In this section, we establish a result simpler that Theorems 1 and   2. We prove that for
$\alpha_0$ small, there is no solution $u(t)$, blowing up in some finite time $T>0$ such that
$$\exists C_2=C_2(u_0),\hbox{ for all } t\in [0,T),\quad 
\int u_x^2(t)\le {C_2\over (T-t)^{2/3}}$$
(note that there is no {\it a priori} control on $C_2$).
Equivalently, for all solutions blowing up in finite time $T>0$, we have
\be\label{weakresult}\mathop{{\rm lim~sup}}_{t\uparrow T} ~ (T-t)^{2/3} \int u_x^2(t) =+\infty.
\ee
The proof is somewhat pedagogical and is based on some geometrical and scaling considerations to
reduce the initial problem to an asymptotic problem; see \S 3.1. The fact that the asymptotic problem has no solution
is then proved in \S 3.2.  This proof of nonexistence, together with the Liouville theorem in \cite{MM2}
will be the base of the proof of Theorem 1.
Nevertheless, we point out that the proof of (\ref{weakresult}) does not rely on the Liouville theorem.

\demo{{\rm 3.1.} Reduction to a solution blowing up as a Dirac
mass in $L^2$}
We prove (\ref{weakresult}) by contradiction.
Assume that  $u(t)$ blows up in finite time $T>0$ and satisfies
$$\hbox{for all } t\in [0,T), \quad  \int u_x^2(t)\le {C_2(u_0)\over (T-t)^{2/3}}.$$
From this assumption, we construct another solution with the same properties as $u(t)$ and blowing up as a Dirac
mass in
$L^2$. This is obtained as a
recurrent object as $t\uparrow T$.\enddemo

\proclaimtitle{Reduction to a solution blowing up as a Dirac mass\break in
$L^2$}
\specialnumber{3}\proclaim{Proposition}\label{REDUCTIONDIRAC}  
For $\alpha_0>0$  small enough{\rm ,} under the preceding assumptions{\rm ,} there exists a solution 
$\widetilde u(t)$ defined on $[0,1)$ such that
\smallbreak
{\rm (i)} $\int \widetilde u^2(t)\le \int Q^2+\alpha_0${\rm ,}
\smallbreak
{\rm (ii)} $E(\widetilde u)\le 0${\rm ,}
\smallbreak
{\rm (iii)} $\widetilde u(t)$ blows up at $t=1${\rm ,} and  for all $t\in [0,1)${\rm ,} 
$|\widetilde u_x(t)|_{L^2}\le {C\over (1-t)^{1/3}}${\rm ,}
\smallbreak
{\rm (iv)} $\widetilde u^2(t,x+\widetilde x(t))\rightharpoonup \left(\int \widetilde u^2(0)\right) \delta_{x=0}
,$ as $t\uparrow 1${\rm .}
\endproclaim

\demo{{P}roof}
Recall that from the assumption on $u(t)$, the introduction and (\ref{17bis}),   we have $\hbox{for all } t\in [0,T),$
$$ \quad  {C_1\over (T-t)^{2/3}}\le \int u_x^2(t)\le {C_2\over (T-t)^{2/3}} \enspace \hbox{ and }\enspace 
{C'_1\over (T-t)^{1/3}}\le
{1\over \l(t)} \le {C'_2\over (T-t)^{1/3}}.\quad$$

Let $(t_n)$ be an increasing sequence of $[0,T)$ such that $t_n\goto T$ as $n\goto +\infty$.
For $t\in [1-T\l_n^{-3}, 1)$, $x\in \R$, we set
\be\label{defofun}
u_n(t,x)=\l_n^{1/2} u(\l_n^3 (t-1)+T,\l_n x ),\enspace\hbox{ where }\enspace 
\l_n=(T-t_n)^{1/3}. \qquad 
\ee
For all $n$, $u_n$ is a solution of (\ref{kdv}) on
$[1-T\l_n^{-3},1)$, blowing up at $t=1$, with initial data
\be\label{defofun0}
u_n(0,x)=\l_n^{1/2} u(t_n,\l_n x),
\ee
and such that
\smallbreak
(i) $|u_n|_{L^2}=|u|_{L^2}$, $|u_n(0)|_{H^1}\le C$ ($C$ independent of $n$),
\smallbreak
(ii) $E(u_n)\goto 0,$ as $n\goto +\infty$,
\smallbreak
(iii) $\hbox{for all } t\in [1-T\l_n^{-3}, 1), \quad
{C_1  \over (1-t)^{2/3}} \le \int (u_n)_x^2(t,x)dx \le {C_2  \over  (1-t)^{2/3}}.$
\smallbreak

Indeed, (i) follows by invariance of the $L^2$ norm by scaling, and 
\be\label{seulementla}{C_1\over (1-t)^{2/3}}={C_1 \l_n^2 \over (\l_n^3 (1-t))^{2/3}}
 \le \int (u_n)_x^2(t,x)dx \le {C_2 \l_n^2 \over (\l_n^3 (1-t))^{2/3}}={C_2\over (1-t)^{2/3}},
\ee
applied to $t=0$.
(ii) follows from $E(u_n)=\l^2_n E(u)\goto 0$ as $n\goto +\infty$.
\enddemo

{\it Step} 1.  Definition and first properties of the limit object.
We claim that there exists $\widetilde u(0)\in H^1(\R)$, $\widetilde u(0)\not \equiv 0$, such that 
$0\le \alpha(\widetilde u(0))\le \alpha_0$, $E(\widetilde u(0))\le 0$, and  
there exists a subsequence of $(u_n)$, still denoted by $(u_n)$, satisfying
$$u_n(0)\rightharpoonup \widetilde u(0)\quad \hbox{in $H^1(\R)$}.$$
 
The existence of a subsequence of $(u_n(0))$ which converges to some $\widetilde u(0)\in H^1(\R)$ 
weakly in $H^1(\R)$ is a consequence of the uniform bound (i)  on $(u_n(0))$ in $H^1$.

Now, let us show that $\widetilde u(0)\not \equiv 0$. We first note that from (\ref{corollaire1}) 
we have, for $\alpha_0$ small enough, $|{x_s\over \l}-1|=|\l^2 x_t-1|\le 1/2$, and so
$$\hbox{for all } t\in [0,T),\quad 0\le x_t(t)\le {3\over 2\l^2(t)}\le {C\over (T-t)^{2/3}},$$
and so $x(t)$ approaches a limit as $t\uparrow T$, denoted $x(T)$. By considering $u(t,x+x(T))$, solution of
(\ref{kdv}) instead of $u(t,x)$, we assume that
$$x(T)=0  \hbox{
and so, by integration,  for all } t\in [0,T),\   |x(t)|\le C (T-t)^{1/3}\le C\l(t).$$
Now, by Lemma \ref{MODULATION},
$|Q-\l^{1/2}(t_n) u(t_n, \l(t_n) x +x(t_n))|_{H^1} \le \delta'(\alpha_0),$
and so, by scaling invariance of the $L^2$-norm,
$$\left|u_n(0,x)-\left({\l_n \over \l(t_n)}\right)^{1/2} Q\left({\l_n \over \l(t_n)} x -{x(t_n) \over
\l(t_n)}  \right)\right|_{L^2}
\le \delta'(\alpha_0).$$
Using ${|x(t_n)|\over \l(t_n)}\le C$ and $C'_1\le {\l_n\over \l(t_n)}\le C_2'$, there exists 
$A$ such that    
\bee |u_n(0)|_{L^2(|x|>A)} \le 2\delta'(\alpha_0),
\eee
which implies for $\alpha_0$ small enough that
$$|u_n(0)|_{L^2(|x|<A)}\ge |Q|_{L^2} -{1\over 4} |Q|_{L^2} ={3\over 4} |Q|_{L^2},$$
and thus $|\widetilde u(0)|_{L^2}\ge {3\over 4} |Q|_{L^2},$  by strong convergence $L^2_{\rm loc}$.

Next, note that by the properties of weak convergence,   $0\le \alpha(\widetilde u(0))\break\le \alpha_0$
and  $E(\widetilde u(0))\le 0.$ The second fact is not obvious, and we prove it briefly.

Let us define a function $\rho$ such that
$$0\le \rho\le 1,~ \rho(x)=1 \hbox{ for $|x|\le 1$, } \rho(x)=0 \hbox{ for $|x|\ge 2$, }
\sqrt{\rho},\sqrt{1-\rho} \in {\cal C}^2,$$
and for $k\in \N$, $\rho_k(x)=\rho({x\over k}).$
First, by direct calculations,
$$E(u_n(0))=E(u_n(0)\sqrt{\rho_k})+E(u_n(0)\sqrt{1-\rho_k})+R_{n,k},$$
where $R_{n,k}=-{1\over 8} \int u_n^2(0)(\rho_k')^2({1\over \rho_k}+{1\over 1-\rho_k})
+{1\over 2} \int u_n^6(0)\rho_k (\rho_k-1)$ so that  
$R_{n,k}\goto R_k=-{1\over 8} \int \widetilde u^2(0)(\rho_k')^2({1\over \rho_k}+{1\over 1-\rho_k})
+{1\over 2} \int \widetilde u^6(0)\rho_k (\rho_k-1)$ as $n\goto +\infty$,
by the  Rellich compactness theorem.

Second, observe that, for all $n$, $|u_n(0)-Q|_{L^2}\le {1\over 2} |Q|_{L^2}$,
and so, there exists $k_0>0$ such that for all  $k>k_0$ and for all $n$, we have
$|u_n(0)\sqrt{1-\rho_k}|_{L^2}\le |Q|_{L^2}.$ Therefore, by the Gagliardo-Nirenberg inequality (\ref{gn}),
we have  $E(u_n(0)\sqrt{1-\rho_k})\break\ge 0$.

Passing to the limit as $n\goto +\infty,$ since $E(u_n(0))\goto 0$, we obtain
$$\overline{\lim}_{n\goto \infty} E(u_n(0)\sqrt{\rho_k})+R_k\le 0.$$
But $\overline{\lim}_{n\goto \infty} E(u_n(0)\sqrt{\rho_k})\ge E(\widetilde u(0)\sqrt{\rho_k}),$
so that 
$E(\widetilde u(0)\sqrt{\rho_k})+R_k\le 0.$
Finally, letting $k\goto +\infty$, we obtain $E(\widetilde u(0))\le 0$.

Now, we consider the solution $\widetilde u(t)$ of (\ref{kdv}) with initial data $\widetilde u(0)$,
defined on $[0,T^*),$ for $T^*>0$.
Since $E(\widetilde u(0))\le 0$ and $\alpha(\widetilde u(0))\le \alpha_0$, 
$\widetilde u(t)$ admits a canonical decomposition as in \S 2.1. 
We denote by $\widetilde \l(t)$, $\widetilde x(t)$ and $\widetilde \e(t)$ the 
parameters of the decomposition.

\medbreak

{\it Step} 2.  Properties of the limit solution.
We consider the asymptotic object $\widetilde u(t)$ instead of the original solution $u(t)$ because
$\widetilde u(t)$ satisfies more properties.
Indeed, we prove that $\widetilde u^2(t)$ concentrates as a Dirac mass at the blow-up time.
We  claim the following result:

\proclaimtitle{The asymptotic solution blows up as a Dirac mass in~$L^2$}
\specialnumber{4}\proclaim{Proposition}\label{DIRAC} 
{\rm (a)} $T^*=1$ and  for all $t\in [0,1)${\rm ,}
$\ds{C_1\over (1-t)^{2/3} }\le \int \widetilde u^2_x(t)\le {C_2\over (1-t)^{2/3}}.$
\smallbreak
\noindent{\rm (b)} $\ds \widetilde u^2(t,x+\widetilde x(t))
\rightharpoonup \left(\int \widetilde u^2(0)\right) \delta_{x=0}$
 as $t\uparrow 1${\rm .} 
\endproclaim

\demo{{P}roof  of {\rm (a)}}
  Let us recall the following result of stability of the weak convergence 
for the KdV equation.

\proclaimtitle{\cite{MM2}}
\specialnumber{6} \proclaim{Lemma} \label{stabilityofweak}
For all $0<t_0<{\rm min}(1,T^*),$ 
\bee
&&
\hbox{for all } t\in [0,t_0],\quad   u_n(t)\rightharpoonup \widetilde u(t) \quad \hbox{in $H^1(\R)$}, \\
&& u_n\goto \widetilde u \quad \hbox{in $C([0,t_0],L^2_{\rm loc}(\R)),$}\qquad \l_n \goto \widetilde \l \quad \hbox{in
$C([0,t_0],\R)$.}
\eee
\endproclaim

{\it Proof.}  See Lemma 8 and Appendix D of \cite{MM2}. Note that   growth
in $H^1$ can be seen in $L^2_{\rm loc}$ norm  by variational arguments. \fin

We have from (iii), (\ref{17bis}) and the preceding lemma,
$$\hbox{for all } t\in [0,{\rm min}(1,T^*)),\quad
{C'_1\over (1-t)^{1/3}}\le {1\over \widetilde \l(t)} \le {C'_2\over (1-t)^{1/3}}.$$
Thus, from (\ref{17bis}),
$$\hbox{for all } t\in [0,{\rm min}(1,T^*)),\quad
{C_1\over (1-t)^{2/3}}\le \int \widetilde u_x^2(t) \le {C_2\over (1-t)^{2/3}};$$
in particular $T^*=1$, and (a) is proved.

\medbreak

{\it Proof of} (b). \enspace
(b) will be a consequence of
\be\label{directly}\hbox{for all } A>0,~\hbox{for all } t\in (0,1), 
\int_{|x|<A(1-t)^{1/3}} \widetilde u^2(t,x) dx \ge \int \widetilde u^2(0,x) dx - {C \over A}-
{C\over A^{3}}.\ee

To prove (\ref{directly}),
let us first note that from uniform estimates on $u_n$, we have

\nonumproclaim{Claim} For all $B>0, n>0${\rm ,} and $t_1,t_2\in [1-T\l_n^{-3},1),$
\be\label{encoreclaim} 
\int_{|x|<B} u_n^2(t_2) \ge \int_{|x|<B/2} u_n^2(t_1) - C {|t_2-t_1|^{1/3}\over B}
-C{|t_2-t_1|\over B^3}.\ee
\endproclaim

\demo{{P}roof of claim {\rm (\ref{encoreclaim})}}  Let $\varphi :\R\goto [0,1]$ be a smooth function such that
$\varphi(x)=1$ for $|x|\le 1/2$, and $\varphi(x)=0$ for $|x|\ge 1$. 
For $B>0$, we set
$\varphi_B(x)=\varphi(x/B)$.

 Now,
\bee
{d\over dt} \int u_n^2(t) \varphi_B(x) dx &=&
-3 \int (u_n)_x^2(t) \varphi_B'+\int u_n^2(t) \varphi_B^{(3)}+{5\over 3} \int u_n^6(t)\varphi_B' \\
&=& -{3\over B} \int (u_n)_x^2(t) \varphi'(x/B) +{1\over B^3} \int u_n^2(t) \varphi^{(3)}(x/B)
\\ &&  +\ {5\over 3B} \int u_n^6(t) \varphi'(x/B).
\eee
Therefore, by conservation of the $L^2$ norm and the Gagliardo-Nirenberg inequality,
\bee \left| {d\over dt} \int u^2(t) \varphi_B(x) dx\right|
&\le &  {C\over B} \int u_x^2(t)+{C\over B^3} +{C\over B} \int u^6(t)|\varphi'(x/B)| \\
&\le & {C\over B^3} + {C\over B} {1\over (1-t)^{2/3}}.
\eee
By integration, we obtain
\bee \int u^2(t_1)\varphi_B -\int u^2(t_2) \varphi_B
&\le&  {C|t_2-t_1| \over B^3} + {C\over B} \left| \int_{t_1}^{t_2} {dt\over (1-t)^{2/3}}  \right|
\\ &\le & {C|t_2-t_1| \over B^3} + {C\over B} \left| \int_{t_1}^{t_2} {dt\over (t_2-t)^{2/3}} \right|
\\ &\le &  C {|t_2-t_1| \over B^3} +C {|t_2-t_1|^{1/3} \over B} .
\eee
Then claim (\ref{encoreclaim}) follows from
$$\int u^2(t_1)\varphi_B \ge \int_{|x|<B/2} u^2(t_1) \quad {\rm and}\quad
\int u^2(t_2)\varphi_B \le \int_{|x|<B} u^2(t_2).$$ 

\smallskip

Now, we prove (\ref{directly}), using the fact that
$\widetilde u$ is recurrent in $u$.
Let $A>0$, and $t\in [0,1)$.
Recall that $\l_n=(1-t_n)^{1/3}$, and
$u_n(t,x)=\l_n^{1/2} u(\l_n^3(t-1)+T,\l_n x).$
When $m\in\N$, we have
$
u(t,x)={1\over \l_{m+n}^{1/2}} u_{m+n}({t-T\over \l_{m+n}^3}+1,{x\over \l_{m+n}}),
$
and so
$$u_n(t,x)={\l_n^{1/2}\over \l_{m+n}^{1/2}} u_{m+n} \left(
{\l_n^3\over \l_{m+n}^3}(t-1)+1,{\l_n \over \l_{m+n}}x \right).$$
By change of variable $y={\l_n\over \l_{m+n}} x$,  
\bee \int_{|x|<A(1-t)^{1/3}} u_n^2(t,x)dx&=&\int_{|y|< {\l_n \over \l_{m+n}} A(1-t)^{1/3}}
u_{m+n}^2\left({\l_n^3\over \l_{m+n}^3}(t-1)+1,y\right)dy.
\eee
Apply claim (\ref{encoreclaim}) to the solution $u_{m+n}$ with $t_2={\l_n^3\over \l_{m+n}^3}(t-1)+1$, $t_1=0$
and $B={\l_n\over \l_{m+n}}A(1-t)^{1/3}$:
\bee &&\hskip-48pt \int_{|y|< {\l_n \over \l_{m+n}} A(1-t)^{1/3}}
u_{m+n}^2\left({\l_n^3\over \l_{m+n}^3}(t-1)+1\right) \\
&&\qquad \ge \int_{|y|<{1\over 2} {\l_n \over \l_{m+n}}A(1-t)^{1/3}}
u_{m+n}^2(0) \\
&&\qquad\quad   -\ C {\left|{\l_n^3\over \l_{m+n}^3}(t-1)+1\right|^{1/3} \over {\l_n \over \l_{m+n}} A(1-t)^{1/3}}-    C
{\left|{\l_n^3\over
\l_{m+n}^3}(t-1)+1\right| \over \left({\l_n \over \l_{m+n}} A(1-t)^{1/3}\right)^3}.
\eee

Now, let 

$\bullet$ $X_0>0$ be such that
$\ds \int_{|x|>X_0} \widetilde u^2(0) \le {1\over A},$

$\bullet$ $n>0$ (depending on $A$, $t$ and $X_0$) be such that
$$\int_{|x|<A(1-t)^{1/3}} \widetilde u^2(t)\ge \int_{|x|<A(1-t)^{1/3}} u_n^2(t)-{1\over A},$$
and 
$$\hbox{for all } k\ge 0,\quad \int_{|x|<X_0} u_{k+n}^2(0) \ge \int_{|x|<X_0} \widetilde u^2(0) -{1\over A},$$

$\bullet$ $m$ be such that
$${1\over 2} {\l_n\over \l_{m+n}} A(1-t)^{1/3} >X_0 \quad {\rm and} \quad
{\l^3_{m+n}\over \l^3_{n}} < {1-t}.$$
Then, we obtain
$${\left|{\l_n^3\over \l_{m+n}^3}(t-1)+1 \right| \over \left({\l_n \over \l_{m+n}} A(1-t)^{1/3}\right)^3}
= {(1-t) {\l_n^3 \over \l_{m+n}^3}-1 \over \left(  {\l_n\over \l_{m+n}} A(1-t)^{1/3}\right)^3}
\le{1\over A^3}.$$
Therefore,
\bee 
\int_{|x|<A(1-t)^{1/3}} \widetilde u^2(t)&\ge& \int_{|x|<A(1-t)^{1/3}} u_n^2(t) -{1\over A}\\ &\ge&
\int_{|x|<X_0}
u_{m+n}^2(0) - {C \over A^3} - {C+1\over A}\\
&\ge &  \int \widetilde u^2(0) - {C\over A^3} -{C+3 \over A};
\eee
thus (\ref{directly}) holds.

Now, we finish the proof of (b).
From (\ref{directly}), we have directly
\be\label{break}\widetilde u^2(t,x)\rightharpoonup
\left(\int \widetilde u^2(0)\right) \delta_{x=0}, \quad \hbox{as $t\uparrow 1$}.\ee
From the fact that $|\widetilde u(t)|_{L^2}\le {C\over (1-t)^{1/3}}$, as in Step 1,
$\widetilde x(t)$
has a limit as $t\uparrow 1$, denoted $\widetilde x(1)$. By     (\ref{break}),
we have  $\widetilde x(1)=0$.
Therefore,
$$\widetilde u^2(t,x+\widetilde x(t))\rightharpoonup
\left(\int \widetilde u^2(0)\right) \delta_{x=0}, \quad \hbox{as $t\uparrow 1$}.$$
Thus Propositions \ref{DIRAC} and \ref{REDUCTIONDIRAC} are proved.\fin

In conclusion, for $\alpha_0$ small enough, from Proposition \ref{REDUCTIONDIRAC},
Lemma \ref{firstexpo} and Corollary \ref{secondexpo}, since $\widetilde \l(t)\goto 0$ as $t\uparrow T$,
we obtain a solution $\widetilde u(t)$ on 
$[0,1)$, such that

\smallbreak

(i) $\int \widetilde u^2 \le \int Q^2+\alpha_0,$
\smallbreak
(ii) $E(\widetilde u)\le 0$,
\smallbreak
(iii) there exists  $s_1<s_2$ satisfying
\begin{eqnarray*}
 \hbox{for all } s\ge s_1,\enspace \widetilde \l(s)&\le& \widetilde \l(s_1),\enspace \widetilde
\l(s_2)={\widetilde
\l(s_1)\over 1.1},\\
 \hbox{for all } s\in [s_1,s_2],\enspace {\widetilde \l(s_1)\over 1.1}&\le& \widetilde \l(s)\le \widetilde \l(s_1),
\end{eqnarray*}
and
$$ \hbox{for all } y<0,\hbox{ for all }  s\in [s_1,s_2],\enspace 
| \widetilde \e(s,y)| \le  C_1 \alpha_0^{1/4} e^{-C_2 |y|}.$$
 
In the next subsection, we prove the nonexistence of such $\widetilde u(t)$, and 
(\ref{weakresult}) follows.

\demo{{\rm 3.2.} Nonexistence of a focusing solution with exponential decay on the left}
In this section, we prove the following proposition, which is crucial for the proof of the 
stability of $Q$ as a blow-up profile. \enddemo

\proclaimtitle{Nonexistence of a focusing solution with exponential\break
decay}
\specialnumber{5}\proclaim{Proposition}\label{NONEXISTENCE} 
There exists $\alpha_I>0$ such that there exists no solution $u(t)$ of {\rm (\ref{kdv})}
satisfying
\smallbreak
{\rm (i)} $\int u^2\le \int Q^2 + \alpha_I${\rm ,}

\smallbreak
{\rm (ii)} $E(u)\le 0${\rm ,}

\smallbreak
{\rm (iii)} There exist $s_1<s_2$ such that
\begin{eqnarray}\label{defs1s2}
&&\hbox{for all } s>s_1,~ \l(s)\le \l(s_1),\enspace 
\l(s_2)={\l(s_1)\over 1.1},\\
 &&\hbox{for all } s\in [s_1,s_2],~ {\l(s_1)\over 1.1} \le \l(s) \le \l(s_1),
\nonumber\end{eqnarray}

\smallbreak
{\rm (iv)}  $\e(s)$ is such that
\be\label{expofinal}
\hbox{for all } y<0, \hbox{ for all } s\in [s_1,s_2],\quad 
|\e(s,y)|\le C_1 \alpha_0^{1/4} e^{-C_2 |y|}.
\ee
\endproclaim

{\it Proof}.
The proof proceeds in three steps. In the first step, we give the main argument assuming
that two fundamental inequalities hold. In steps 2 and   3, we prove these two inequalities.  We work on the time
interval
$[s_1,s_2]$.

\demo{Step 1. {\it The main argument}}
The problem in this proof is to understand the size of different quantities such as
$$\int_{s_1}^{s_2} \int \e^2,\quad \int_{s_1}^{s_2} \int \e_y^2,\quad
\int_{s_1}^{s_2} \int \e^2 e^{-{|y|\over 2}}$$
(and also at some point $\int_{s_1}^{s_2} \int \e Q.$) Two difficulties are  that we have no control
of the size of the interval $[s_1,s_2]$ (the ``doubling time") in terms of $\alpha_0$, and that the
terms we integrate in time are oscillatory integrals.
\enddemo

We claim that there exists $\alpha_7>0$ such that if $\alpha_0<\alpha_7$ then the following two inequalities hold:

\medbreak
$\bullet$ There exists $C_I>0$ (independent of $\alpha_0$) such that 
\be\label{etoile1}
C_I \int_{s_1}^{s_2} \int \e^2 e^{-{|y|\over 2}} \ge 1+\int_{s_1}^{s_2} \int \e_y^2
+|E_0| \int_{s_1}^{s_2} \l^2.
\ee
This is a consequence of a dispersion relation in $L^1$ in terms of $J(s)$ 
(Proposition~\ref{FORMULEJ}) and exponential decay. See step 2 for the proof.

\medbreak

$\bullet$ Let $A_0>2$ be as in Proposition~2. 
There exists $C_{II}>0$ (independent of $\alpha_0$) such that 
\be\label{etoile2}
 \int_{s_1}^{s_2} \int \e^2 e^{-{|y|\over A_0}} \le C_{II}\,\alpha_0^{1/2} \left[ 1+\int_{s_1}^{s_2} \int \e_y^2
+|E_0| \int_{s_1}^{s_2} \l^2\right].
\ee
This is a consequence of the local Virial type relation (Proposition~2), i.e.\ $L^2$ dispersion. See step 3 for the proof.
\medskip

Note that since $A_0> 2$,  the contradiction is obvious for $\alpha_0<\left({1 \over C_I C_{II}}\right)^{2}$.
Thus we need only prove (\ref{etoile1}) and (\ref{etoile2}).

\demo{Step 2.  {\it Proof of} (\ref{etoile1})}
Recall that we have the following definition of~$J(s)$:
$$J(s)=\int \e(s,y)\int_y^{+\infty} \left({Q\over 2}+z Q_z\right) - {1\over 4} \left(\int Q\right)^{2}.$$
Also,  let
$$\widetilde J(s)=\int \e(s,y)\int_y^{+\infty} \left({Q\over 2}+z Q_z\right).$$
By the exponential decay on the left (\ref{expofinal}), $J(s)$ is well-defined for all $s\ge 0$, and 
it follows from (\ref{formuleJ}) that 
$$J_s+{\l_s\over 2\l} J+2\int \e Q \le C \int \e^2 e^{-{|y|\over 2}}.$$
By (\ref{energie}), for $\alpha_0$ small enough ($E_0\le 0$),
$$2\int \e Q \ge 2 \l^2 |E_0| +{1\over 2} \int \e_y^2-C \int \e^2 e^{-|y|}.$$
Therefore,
$$J_s+{\l_s\over 2\l} J+{1\over 2} \int \e_y^2 +2 \l^2 |E_0| \le C \int \e^2 e^{-{|y|\over 2}}.$$
Multiplying this relation by $\sqrt{\l(s)\over \l(s_1)}$ we find
$$\left(\sqrt{\l(s)\over \l(s_1)} J \right)_s
+{1\over 2} \sqrt{\l(s)\over \l(s_1)} \int \e_y^2 +2 \sqrt{\l(s)\over \l(s_1)} \l^2 |E_0|
\le C \sqrt{\l(s)\over \l(s_1)} \int \e^2 e^{-{|y|\over 2}}.$$
Since for $s\in [s_1,s_2]$, $1\ge \sqrt{\l(s)\over \l(s_1)}\ge {1\over \sqrt{1.1}}\ge {1\over 2},$
we obtain
$$\left(\sqrt{\l(s)\over \l(s_1)} J \right)_s
+{1\over 4}  \int \e_y^2 + \l^2 |E_0|
\le C  \int \e^2 e^{-{|y|\over2}}.$$
After integration between $s_1$ and $s_2$, we obtain
$$\left[\sqrt{\l(s_2)\over \l(s_1)} J (s_2)- J(s_1)\right]
+{1\over 4} \int_{s_1}^{s_2} \int \e_y^2 + |E_0| \int_{s_1}^{s_2} \l^2
\le C \int_{s_1}^{s_2} \int \e^2 e^{-{|y|\over 2}}.$$

Next, we have
$$\sqrt{\l(s_2)\over \l(s_1)} J (s_2)- J(s_1)
={1\over 4} \left(1-\sqrt{1\over 1.1}\right) \left(\int Q\right)^2
+\sqrt{1\over 1.1} \widetilde J(s_2)-\widetilde J(s_1).$$
Note that from the exponential control on the left of $\e(s,y)$,
and the decay properties of the functions $Q$ and $Q_y$,
for all  $s\in [s_1,s_2]$,
$$\left|\widetilde J(s)\right|\le C \alpha_0^{1/4} \left| \int_y^{+\infty} {Q\over 2}+zQ_z\right|_{L^\infty}
\int_{y<0} e^{C_2 y} +C |\e(s)|_{L^\infty} \int_{y>0} e^{-C y} \le C \alpha_0^{1/4}.$$
Therefore, for $\alpha_0$ small enough,
\begin{eqnarray*}
\sqrt{\l(s_2)\over \l(s_1)} J (s_2)- J(s_1)
&\ge& {1\over 4}\left(1-\sqrt{1\over 1.1}\right) \left(\int Q\right)^2 - 2C \alpha_0^{1/4}
\\
&\ge&{1\over 8}\left(1-\sqrt{1\over 1.1}\right) \left(\int Q\right)^2.
\end{eqnarray*}

Thus, there exists $C_I>0$ such that
$$C_I \int_{s_1}^{s_2} \int \e^2 e^{-{|y|\over 2}} \ge 1+\int_{s_1}^{s_2} \int \e_y^2
+|E_0| \int_{s_1}^{s_2} \l^2.
$$
\enddemo

\demo{Step 3.  {\it Proof of} (\ref{etoile2})}
By Proposition~2, and the orthogonality conditions imposed on $\e(s)$
(see Lemma \ref{MODULATION}), 
$$\left (\int \Psi_{A_0} \e^2\right)_s
\le -\delta_0  \int (\e^2 +\e_y^2) e^{-{|y|\over A_0}}
+{1 \over \delta_0} \left(\int \e Q\right)^2.$$
Integrating between $s_1$ and $s_2$, we obtain:
$$\int \Psi_{A_0} \e^2(s_2) - \int \Psi_{A_0} \e^2(s_1)
\le -{\delta_0} \int_{s_1}^{s_2} \int (\e^2 +\e_y^2) e^{-{|y|\over A_0}}
+{1 \over \delta_0} \int_{s_1}^{s_2} \left(\int \e Q\right)^2.$$

On the one hand, we have
$$\left|  \int \Psi_{A_0} \e^2(s_2) - \int \Psi_{A_0} \e^2(s_1)\right|
\le C A_0 |\e(s)|^2_{L^2}\le C' \alpha_0,$$
since $|\Psi_{A_0}(y)|\le C A_0$ and $\int \e^2\le C\alpha_0$.

On the other hand, by $\int \e^2\le C\alpha_0$, 
$$\int_{s_1}^{s_2} \left( \int \e Q\right)^2 \le \sup_{(s_1,s_2)} \left(\int |\e Q|\right)
\int_{s_1}^{s_2} \left|\int \e Q\right|\le C \alpha_0^{1/2} \int_{s_1}^{s_2} \left|\int \e Q\right|.$$
By (\ref{energie}), 
$$\left| \int \e Q\right| \le {\l^2} |E_0| +\int \e_y^2 +C \int \e^2 e^{-|y|}.
$$
Therefore,
\bee
\int_{s_1}^{s_2} \left( \int \e Q\right)^2 
&\le& C\alpha_0^{1/2} \int_{s_1}^{s_2} \left(\l^2 |E_0| +\int \e_y^2 + \int \e^2 e^{-|y|}\right)
\\ &\le&  C\alpha_0^{1/2} \int_{s_1}^{s_2} \left(\l^2 |E_0|+\int \e_y^2\right) +{\delta_0^2\over 2}
\int_{s_1}^{s_2} \int \e^2 e^{-{|y|\over A_0}},
\eee
for $\alpha_0$ such that $C \alpha_0^{1/2}\le \delta_0^2/2$.
\pagebreak

Finally, we obtain, for $\alpha_0$ small,
$$\int_{s_1}^{s_2} \int(\e^2+\e_y^2) e^{-{|y|\over A_0}}
\le C_{II} \left[ {\alpha_0}+\alpha_0^{1/2} \int_{s_1}^{s_2}\left(\l^2 |E_0| +\int \e_y^2 \right)
\right];$$
thus the proof of (\ref{etoile2}) is complete, and Proposition \ref{NONEXISTENCE} is proved.\hfill\qed
\enddemo

\section{Lower bounds on the blow-up rate and\\ stability of $Q$ as blow-up profile}
\advance\eqcount by 43

In this section, we prove the two main results. First, we show that Theorem 2 follows
easily from Theorem 1. The rest of the section is then devoted to the proof of Theorem 1.

\demo{{\rm 4.1.} Theorem {\rm 1} implies Theorem {\rm 2}}
The proof of the lower bounds on the blow-up rate follows directly from the fact that in the
rescaled variables, the solution behaves like $Q$, and the fact that the solution 
$v(t,x)=Q(x-t)$ satisfies ${d \l_v \over dt}=0$.

Consider a solution $u(t)$ blowing up in a finite time $T>0$,
with\break $\int u^2<\int Q^2+\alpha_0$, where $\alpha_0$ is defined as  in Theorem 1.
From Theorem 1, we have
$$\e(s) \rightharpoonup 0,\quad \hbox{in $H^1(\R)$, as $s\goto +\infty$}.$$
This implies, by the compact embedding from $H^1(\R)$ to $L^2_{\rm loc}(\R)$, that 
$$\e(s) e^{-{|y|\over \sqrt{2}}}\goto 0,\quad \hbox{in $L^2(\R)$ strong, as $s\goto +\infty$}.$$

From Lemma \ref{PROPERTIESDECOMPOSITION} (iii),  we have
$ \left|{\l_s\over \l}\right| \le C \left(\int \e^2(s) e^{-{|y|\over 2}}\right)^{1/2},$
and so 
${\l_s(s)\over \l(s)}\goto 0,$ as $s\goto +\infty$.
Since ${ds\over dt}=1/\l^3$, this is equivalent to
$$(\l^3(t))_t=3 \l^2(t) \l_t(t)\goto 0,\quad \hbox{as $t\uparrow T$.}$$
By integration, it follows that 
${\l^3(t)\over (T-t)}\goto 0,$ as $t\uparrow T$,
and then (\ref{17bis}) yields Theorem~2.\enddemo

4.2. {\it Proof of the stability of the blow-up profile}.
In this subsection, we prove Theorem 1.
We consider a solution $u(t)$ of (\ref{kdv}) blowing up in  finite 
or infinite time $T>0$. We assume that 
$\int u^2 =\int Q^2 +\alpha_0$ where $\alpha_0$ is  small so that we have 
the decomposition described in Section~2.
Since $\l^{1/2}(s) u(s,\l(s)y+x(s))=Q(y)+\e(s,y)$, the result is implied by
\be\label{theresult}\e(s)\rightharpoonup 0 \quad \hbox{in $H^1(\R)$, as $s\goto +\infty$.}
\ee

The proof of (\ref{theresult}) is divided in two parts. 
\pagebreak

In part A, we prove (\ref{theresult})
for a subsequence of time $(t_n)$ related to a monotonicity property of the oscillations of
$|u_x(t)|_{L^2}$ or equivalently $\l(t)$, namely
for all  $t\ge t_n$, $\l(t)\le \l(t_n)$, $\l(t_{n+1})={\l(t_n)\over 1.1}.$
Under this control, we are able to 
construct a limit object $\widetilde u(t)$
with uniform exponential decay for $y<0$; this property is one of the keys to obtain the
contradiction with $\widetilde u(t)$ as in Section~3. 
However, unlike Section~3, we do not have    control of the growth of the gradient
$|u_x(t)|_{L^2}\le {C\over (T-t)^{1/3}}.$ Thus,  
two possibilities will have to be considered for $\widetilde u(t)$:

\smallbreak

(i) Either there exists $t>0$ such that $\widetilde \l(t)={1\over 1.1}$. Thus the solution grows in $H^1$ for $t>0$
and then we are reduced to the case of
\S 3.2;

\smallbreak

(ii) Or for all $t>0$, ${1\over 1.1}< \widetilde \l(t)\le 1$ which means that
the solution stays bounded, and the Liouville theorem of
\cite{MM2}  giving a classification  of such solutions gives a  contradiction. 
More precisely, we recall the following result of  asymptotic stability of $Q$.

\proclaimtitle{Asymptotic stability of $Q$ \cite{MM2}}
\specialnumber{6}\proclaim{Proposition}\label{LIOU}
There exists $\alpha_{II}>0$ such that if a  solution $u(t)$ of {\rm (\ref{kdv})}  defined for $t\ge 0$ is such that
\smallbreak
{\rm (i)} $\int u^2<\int Q^2+\alpha_{II},$

\smallbreak
{\rm (ii)} $E(u)\le 0${\rm ,}

\smallbreak
{\rm (iii)} $\hbox{for all } t\ge 0$, $\ds {1\over 2} \le \l(t)\le 1${\rm ,}
\smallbreak
\noindent then $\e(t)\rightharpoonup 0$ in $H^1(\R)$, as $t\goto +\infty${\rm .}
\endproclaim

This is a direct consequence of Theorem 2 in \cite{MM2}, and Lemma \ref{MODULATION}. Note that the
result is strongly related to a classification of bounded solutions close to $Q$, characterizing the
soliton $Q(x-t)$.

In both (i) and (ii), $L^1$ dispersion (Proposition \ref{FORMULEJ}) and $L^2$ dispersion (Proposition~2) play a crucial role to obtain a contradiction.

In part B, we 
extend the result to the whole sequence $t\goto T$, using arguments in $L^2$. More precisely,
we use almost-monotonicity of the mass to the left of the solution when $t$ increases, and to the right of the
solution when $t$ decreases (using the argument backwards, and the invariance of equation (\ref{kdv}):
if $u(t,x)$ is solution then $u(-t,-x)$ is solution). Therefore, the $L^2$ mass of 
$Q+\e(t)$ on compact sets in space is controlled for all time, and we conclude by using the 
variational characterization of $Q$.

\demo{{\rm A.} Convergence to $Q$ as $t_n\goto T$}
We define an increasing sequence $t_n\goto T$
such that 
\be\label{99}\l(t_n)={1\over (1.1)^n},\quad \hbox{for all } t\in [t_n,T), ~ \l(t)\le \l(t_n).
\ee
Such a sequence $(t_n)$ allows us to obtain the property of exponential decay on the left for the limit
object. Indeed, 
such a property is known to hold when there is a control  of   type (\ref{99}) on the solution
(see \S 2.3).

The objective of this subsection is to prove the following proposition.

\proclaimtitle{Convergence as $t_n\goto T$}
\specialnumber{7}\proclaim{Proposition}\label{propclaimA} 
There exists $\alpha_{III}>0$ such that if $\int u^2\le \int Q^2+\alpha_{III}${\rm ,} then we have
\be\label{claimA}
\e(t_n)\rightharpoonup 0\quad \hbox{
in $H^1(\R)$ weak{\rm ,} as $n\goto +\infty$.}
\ee
\endproclaim

{\it Proof of Proposition {\rm \ref{propclaimA}}}.

\demo{Step 1. {\it Introduction of the limit dynamics}}
For $\overline x\in \R$ and $\overline t\in [-{t_n\over \l^3(t_n)},\break {T-t_n \over \l^3(t_n)}),$
we define
$$u_n(\overline t,\overline x)=\l^{1/2}(t_n) u(t_n+\l^3(t_n)\overline t,
\l(t_n) \overline x+x(t_n)).$$
By scaling invariance of (\ref{kdv}), the function $u_n(\overline t,\overline x)$ is a solution of (\ref{kdv}) with 
blow-up time $T_n={T-t_n \over \l^3(t_n)}$ ($T_n=+\infty$ if $T=+\infty$). Note that
${T-t_n\over \l^3(t_n)}\goto +\infty$ as $n\goto +\infty$ is possible.

We check easily that the sequence $(u_n)$ satisfies the following properties:

\smallbreak
(i) $\int u_n^2 =\int Q^2+\alpha_0$,
\smallbreak
(ii) $\lim_{n\goto +\infty} E(u_n)=0$ (since $E(u_n)=\l^2(t_n) E(u)$),
\smallbreak
(iii)  $\l_n(0)=1$, $x_n(0)=0$, $\hbox{for all } \overline t\in [0,T_n),$ $\l_n(\overline t)\le 1$,
\smallbreak
(iv) $|u_n(0)|_{H^1}\le C$,
\smallbreak \noindent  where $\l_n(\overline t)$ and $x_n(\overline t)$ are the geometrical parameters associated
to $u_n$. Note that $\e(t_n)=\e_n(0),$ where $\e_n$ is associated to $u_n$ as in \S 2.1.
\enddemo

The proof of the proposition is by contradiction. 
Assume that
up to the extraction of a subsequence (still denoted $\e_n$), there exists 
$\widetilde \e(0)\not \equiv 0$, $\widetilde \e(0) \in H^1(\R)$ such that 
$$\e_n(0)\rightharpoonup \widetilde \e(0)\quad \hbox{in $H^1(\R)$, as $n\goto +\infty$.}$$

We denote by $\widetilde u(t)$ the solution of (\ref{kdv}) with initial data 
$\widetilde u(0)=Q+\widetilde \e(0)$, and  $\widetilde T>0$ its maximal time of existence.
We have easily the following properties of $\widetilde u$.

\proclaimtitle{First properties of the limit object $\widetilde u$}
\specialnumber{7} \proclaim{Lemma} \label{proputilde}
We have
\smallbreak
{\rm (i)} $\int Q^2\le \int \widetilde u^2\le \int Q^2+\alpha_0${\rm ,}

\smallbreak
{\rm (ii)} $E(\widetilde u)\le 0${\rm ,}

\smallbreak
{\rm (iii)} $\widetilde \l(0)=1${\rm ,} $\widetilde x(0)=0${\rm ,}  for all $t\in [0,\widetilde T),$ $ \widetilde \l(t)\le 1${\rm .}
\endproclaim

\demo{{P}roof} 
(i) and 
(ii) follow from standard variational arguments given in \S 3.1.
\smallbreak
(iii) Since $E(\widetilde u)\le 0$, we can consider $\widetilde \e$, $\widetilde \l$ and $\widetilde x$ 
associated to the
decomposition of $\widetilde u$. Since $\int xQ_x \e_n(0)=\int  x \left({\ts {Q\over 2}}+xQ_x\right)\e_n(0)=0,$
by weak convergence, this is still true for $\widetilde \e(0)$. Since
$\widetilde u(0)=Q+\widetilde \e(0)$,  
$\widetilde \l(0)=1$ and $\widetilde x(0)=0$.

We recall from Lemma \ref{stabilityofweak} that for any $t_0\in [0,{\rm min}(\widetilde T,\underline{\lim} {T_n})),$
\bee &&
\hbox{for all } t\in [0,t_0],\quad  u_n(t)\rightharpoonup \widetilde u(t) \quad \hbox{in $H^1(\R)$,}
\\ &&  u_n\goto \widetilde u \quad
\hbox{in $C([0,t_0],L^2_{\rm loc}(\R)),$} \qquad  \l_n \goto \widetilde \l \quad \hbox{in $C([0,t_0],\R)$ as $n\to +\infty$}.
\eee

We claim that $\widetilde T\le \underline{\lim}{T_n}$. By contradiction, $\widetilde T> \underline{\lim}{T_n}$ and
the preceding property imply that there exists $c>0$
such that for all  $t_0\in [0,\underline{\lim}{T_n}),$  $$\exists n(t_0), \hbox{ for all } n\ge n(t_0),\l_n(t_0)\ge c.$$
By (\ref{17bis}) it follows that $|u_{nx}(t_0)|_{L^2}\le C$; thus by the well-posedness of the 
Cauchy problem in $H^1$, $u_n(t)$ is defined on $[t_0,t_0+\tau]$, where $\tau>0$ is independent
of~$t_0$, and in particular $\underline{\lim}{T_n}\ge t_0+\tau.$ This is a contradiction.

Thus, we obtain
for all $t\in [0,\widetilde T),$ $\widetilde \l(t)\le 1$ and
Lemma \ref{proputilde} is proved. \hfill\qed\enddemo

Now we define
$\tau={\rm inf}\{t>0, \widetilde \l(t)={1\over 1.1}\},$
possibly $\tau=+\infty$.
It follows from the definition of $\tau$ that
$$\hbox{for all } t\in [0,\tau),\quad {1\over 1.1}\le \widetilde \l(t)\le 1.$$

Using  crucially the control of $\widetilde \l(t)$, we have the following lemma, proved in the last step.

 \proclaimtitle{Exponential decay of $\widetilde \e$ on the left on
$[0,\tau)$}
\specialnumber{8} \proclaim{Lemma}\label{EXPOEPS} There exists $C>0$ such that for $\alpha_0$ small enough{\rm ,} 
\vglue4pt
\hfill ${\displaystyle \hbox{for all } y<0,\hbox{ for all } t\in [0,\tau),\quad
|\widetilde \e(t,y)| \le C \alpha_0^{1/4}  e^{-{|y|\over 28}}.}$\hfill
\endproclaim

{\it Step} 2.  The main argument.
The idea is to use two crucial dynamical arguments which say that 
a solution different from $Q$ has to disperse some mass on the left. Therefore, we have a result of
nonexistence of $\widetilde u$, which is a contradiction.
Consider two cases.
 
\demo{{\rm (a)} Focusing regime under exponential decay on the left: $\tau <+\infty$}
The contradiction follows from Proposition \ref{NONEXISTENCE}, for $\alpha_0$ small.
Indeed, from Lemma \ref{proputilde} and Lemma \ref{EXPOEPS}, we obtain a solution $\widetilde
u(t)$ satisfying

\smallbreak

(i) $\int \widetilde u^2(t)\le \int Q^2 +\alpha_0$,
\smallbreak

(ii) $E(\widetilde u)\le 0$,
\pagebreak

(iii) $\widetilde \l(0)=1$ and
$\hbox{for all } t\in [0,\widetilde T), ~\widetilde \l(t)\le 1; \quad \widetilde \l(\tau)={1\over 1.1},  \hbox{ for all } t\in
[0,\tau),~{1\over 1.1}\le
\widetilde \l(t)\le 1$,
\smallbreak

(iv) $\hbox{for all } y<0,\hbox{ for all } t\in [0,\tau),\quad
|\widetilde \e(t,y)| \le  C \alpha_0^{1/4}  e^{-{|y|\over 28}}.$
\enddemo

\demo{{\rm (b)} Regular regime under $L^2$ compactness of the solution: $\tau =+\infty$}
We have a  solution
$\widetilde u(t)$, defined for  $t\ge 0$  such that 
$$\hbox{for all } t\ge 0, \quad {1\over 1.1}\le \widetilde \l(t)\le 1.$$
Using Proposition \ref{LIOU} and the compact embedding of $H^1(\R)$ in $L^2_{\rm loc}(\R)$, we have first:
$$\widetilde \e(t)\goto 0 \quad \hbox{in $L^2_{\rm loc}(\R)$, as $t\goto +\infty$}.$$
\enddemo
Now, using the exponential decay and a property of equation (\ref{kdv}), we have

\nonumproclaim{Claim}  
For all $\delta>0${\rm ,} there exist $R(\delta)>0,\hbox{ for all } t\ge 0$, 
$$\int_{|y|>R(\delta)} \widetilde \e^2(t,y)dy \le \delta.$$
\endproclaim

{\it Proof}.
On the one hand, by the uniform exponential decay of $\widetilde \e$ on the left given by Lemma \ref{EXPOEPS},
we have  the $L^2$ compactness of $\widetilde \e(s)$  on
the left (i.e.\ for $y<0$).

On the other hand, it follows by using backwards the almost monotonicity of the mass that $\widetilde \e$ is 
$L^2$ compact on the right. 
Indeed, assume for the sake of contradiction that there exists $\delta>0$ and a sequence $\widetilde t_n$  be such that 
$\int_{y>n} \widetilde \e^{\,2}(\widetilde t_n, y) dy\ge \delta.$
Since $t\mapsto \widetilde \e(t)$ is continuous in $L^2$, we necessarily have $\widetilde t_n\goto +\infty$.
Let $x_0$ be such that $C_0 e^{-{x_0\over 3}}<\delta/2$.
Since, by the invariance of (\ref{kdv}),
 $\widetilde u(\widetilde t_n-\overline t,-x)$ is a solution, for $\overline t\in
[0,\widetilde t_n]$, by using Lemma \ref{monotonicite}, we have 
${\cal I}_{x_0}(\widetilde t_n)\ge {\cal I}_{x_0}(0)-C_0 e^{-x_0/3}.$ Thus
${\cal I}_{x_0}(\widetilde t_n)\ge \delta-\delta/2=\delta/2$, and for
$n$ large we obtain $\int_{y>n} \widetilde \e^{\,2}(0,y) dy \ge {\delta\over 4}.$ 
But this is a contradiction for $n$ large.
Thus the claim  is proved. 
 \vglue12pt

In conclusion,  $\widetilde \e(t)\goto 0$  in $L^2(\R)$ as $t\goto +\infty$, and 
passing to the limit as $t\goto +\infty$ in the mass conservation
$\int \widetilde u^2(0)=\int \widetilde u^2(t)=\int Q^2+2\int \widetilde \e(t) Q +\int \widetilde \e^{\,2}(t)$, we obtain 
$$\int \widetilde u^2(0)=\int Q^2.$$
Since $E(\widetilde u(0))\le 0$, by the characterization of $Q$ (\ref{charac}), there exist 
$\l_0>0$ and $x_0\in \R$ such that $\widetilde u(0)=Q+\widetilde \e(0)=\l_0^{1/2} Q(\l_0 (x+x_0))$.
Since $|\l_0-1|+|x_0|\le C\sqrt{\alpha_0}$ (from $|\widetilde \e(0)|_{L^2}\le C\sqrt{\alpha_0}$),
 we obtain  by the
orthogonality conditions on $\widetilde \e(0)$ that $\l_0=1$ and $x_0=0$; thus 
 $\widetilde \e(0)\equiv 0$. This is a contradiction.

Therefore,  we have a contradiction in cases (a), (b), and the proposition is proved.

\demo{Step 3. {\it Exponential decay of $\widetilde \e$ on the left on a {\rm ``}\/doubling\/{\rm "}\/ interval of time}}
This step is devoted to the proof of
 Lemma \ref{EXPOEPS}. We claim that Lemma~8 is a consequence of 
\be\label{suru}
\exists C_1>0,\hbox{ for all } x_0<0,\hbox{ for all } t\in [0,\tau), 
\int_{2x_0<x<x_0} \widetilde u^2(t,x+\widetilde x(t))dx \le {C_1} e^{-{|x_0|\over 12}}.
\ee
Indeed,  by summation of (\ref{suru}), properties of $Q$ and control of $\widetilde \l(t)$ on $[0,\tau)$, 
we have directly
$$\hbox{for all } y_0<0,\hbox{ for all } t\in [0,\tau),\quad
\int_{y<y_0} \widetilde \e^2(t,y) dy\le  C e^{-{|y_0|\over 14}},$$
and we obtain Lemma \ref{EXPOEPS} as in the proof of Corollary \ref{secondexpo}.
\enddemo

Now, we prove the exponential decay on $\widetilde u$ (\ref{suru}).
As in the proof of the exponential decay in
\S 3.1, it follows from monotonicity properties on $u_n$ in $L^2_{\rm loc}$ and recurrence of $\widetilde u$
on $u_n$. Here, we do not assume $|u_x(t)|_{L^2}\le {C\over (T-t)^{1/3}},$ and thus $\widetilde u$
does not concentrate as a Dirac mass in $L^2$ at the blow-up time. However, considering a quantity
$\cal J$ which measures the mass lost on the left as $n\goto +\infty$, we still have similar estimates as in \S 3.1.

We prove (\ref{suru}) by contradiction.
Let $a_0$ and $C_0$ be defined as in Lemma~\ref{monotonicite}. Let $a_1>a_0$ and $C_1>0$ be chosen later.
It is sufficient to prove (\ref{suru}) for $x_0<-a_1$.
 Assume that there exists $t_0\in [0,\tau)$, and
$x_0<-a_1$ such that 
$$\int_{2x_0<x<x_0} \widetilde u^2(t_0,x+\widetilde x(t_0)) dx\ge C_1 e^{-{|x_0|\over 12}}.\qquad $$
Since $u_n(t_0,x+x_n(t_0))\goto \widetilde u(t_0,x+\widetilde x(t_0))$ in $L^2_{\rm loc}(\R)$, 
there exists $n_0$ such that 
\be\label{contra}\hbox{for all } n\ge n_0,\quad 
\int_{2x_0<x<x_0} u_n^2(t_0,x+x_n(t_0)) dx \ge {C_1\over 2} e^{-{|x_0|\over 12}}.\ee

For the same choice of function $\psi$ as in section 2.3, for $x_1\in \R$, we define
for $t\in [t_0,T_n)$,
$${\cal J}_{x_1,n}(t)=\int u_n^2(t,x) (1-\psi(x-x_n(t_0)-x_1-{\ts {1\over 4}}(x_n(t)-x_n(t_0)))) dx.$$
Note that this quantity measures the mass on the left since $\psi(x)\goto 0$ as
$x\goto -\infty$, and $\psi(x)\goto 1$ as $x\goto +\infty$.
We claim the following lemma, proved in Appendix B.

\proclaimtitle{Limit of ${\cal J}_{x_1,n}(t)$}
\specialnumber{9} \proclaim{Lemma} \label{suiteJ}
 There exists $\cal J$ such that the following property holds\/{\rm :}
$\exists C>0${\rm ,} $\hbox{for all } \delta_1>0${\rm ,} $\hbox{for all } x_1<-a_1${\rm ,}
$\exists n_1=n_1(\delta_1,x_1)$ such that
\be\label{trois}\hbox{for all } n\ge n_1,\hbox{ for all } t\in [t_0,T_n),\quad 
|{\cal J}_{x_1,n}(t)-{\cal J}| \le \delta_1 + C e^{x_1\over 8}.\qquad
\ee
\endproclaim

\vglue-18pt
{\it Remark}. In fact, $\cal J$ measures the mass lost at the blow-up time at the left of the
soliton in the rescaled variable.
\medbreak

This result provides a contradiction and thus proves (\ref{suru}).
Indeed, on the one hand, for $n\ge n_0$,
we have
\bee {\cal J}_{{x_0},n}(t_0)-{\cal J}_{{2 x_0},n}(t_0) 
& = & \int u_n^2(t_0,x+x_n(t_0)) 
\left(\psi(x-2x_0)-\psi(x-x_0)\right) dx.\eee
Since the function $\psi$ is increasing, and since for $|x_0|$ large enough,
$$\min_{(2x_0,x_0)}(\psi(x-2x_0)-\psi(x-x_0)) \ge {1\over 4},$$
we obtain, by taking $a_1$ large enough and using (\ref{contra}),
\bee {\cal J}_{x_0 ,n}(t_0)-{\cal J}_{{2 x_0},n}(t_0) & \ge & {1\over 4} \int_{2x_0<x<x_0} u_n^2(t_0,x+x_n(t)) dx \ge
 {C_1\over 8} e^{-{|x_0|\over 12}}.\eee

On the other hand, by Lemma \ref{suiteJ}, applied for $\delta_1={C_1\over 32} e^{x_0\over 12}$, 
there exists
$n_1=n_1(x_0)$ such that for all $ n\ge n_1$,
$$
{\cal J}-{C_1\over 32} e^{x_0\over 12}-C e^{{2x_0\over 8}} \le {\cal J}_{2x_0,n}(t_0)
\quad {\rm and} \quad {\cal J}_{x_0,n}(t_0)\le {\cal J}
+{C_1\over 32}
 e^{x_0\over 12} + C e^{x_0\over 8}.$$ 
Therefore,
${\cal J}_{x_0,n}(t_0)-{\cal J}_{{2 x_0},n}(t_0)\le
 {C_1\over 16}
 e^{x_0\over 12} + 2C e^{x_0\over 8}.$
Thus, 
we obtain a contradiction by taking  $C_1>32 C$ (recall that $C$ here is independent of $x_0$)
for $n>n_0,n_1$.
\hfill\qed

\demo{{\rm B.} Convergence to $Q$ for $t\goto T$}
Note that for $(t_n)$ defined in part A,  $t_1<t_2<\ldots<t_n<T,$ and $t_n\goto T$.
We now claim the following proposition, which completes the proof of Theorem 1.
\enddemo

\proclaimtitle{Convergence for $t\goto T$}
\specialnumber{8}\proclaim{Proposition}\label{wholeseq}
There exists $\alpha_{IV}>0$ such that if $\int u^2\le \int Q^2 +\alpha_{IV}${\rm ,} then
$$\e(t)\rightharpoonup 0 \quad \hbox{in $H^1(\R)$, as $t\uparrow T$.}$$
\endproclaim

{\it Proof}.
We argue by contradiction. Suppose that there exists
$\widetilde \e \in H^1(\R)$, $\widetilde \e\not \equiv 0$, and  a subsequence of $(t_n)$ 
denoted $(t_{n'})$, such that for all $n'$ there is $t^1_{n'} \in [t_{n'},t_{n'+1}]$,
such that
$$\e(t^1_{n'})\rightharpoonup \widetilde \e \quad \hbox{in $H^1(\R)$ as $n'\goto +\infty.$}$$

First, we use again the fact that $u_n$ is recurrent in $u$ and the mass properties (Lemma \ref{suiteJ})
to obtain that the $L^2$ mass of $u(t_n)$ as $n\goto +\infty$ has a limit at the left
($x<x(t_n)$) and at the right ($x>x(t_n)$).

Then, we use the almost-monotonicity property twice to extend this result to $u(t^1_{n'})$ 
as $n'\goto +\infty$:
\smallbreak
(i) the first time on $u(t,x)$, from $t_{n'}$ to $t^1_{n'}$ for $x<0$.
\smallbreak
(ii) the second time on $u(-t,-x)$, which is also a solution. $L^2$ estimates at
$t_{n'+1}$ give $L^2$ estimates at $t^1_{n'}$  for $x>0$. Here we use crucially
$\l(t_{n'+1})={\l(t_{n'})\over 1.1},$ to be able to use the  almost-monotonicity property on the interval
$[t_{n'},t_{n'+1}]$.
\smallbreak
We then conclude by energy arguments (variational characterization of
$Q$; see (\ref{charac})).
\vglue12pt

{\it Step} 1.  Limit of the mass  of $\e(t_n)$ as $n\goto +\infty$.
From 
Lemma \ref{suiteJ} and Proposition~\ref{propclaimA}, for ${\cal J}$  defined  as in Lemma \ref{suiteJ},
we claim  that for all $y_1\in
\R$, 
\be\label{masslimit}
\lim_{n\goto +\infty} \int_{y<y_1} \e^2(t_n,y) dy ={\cal J}\quad
\hbox{and}\quad  \lim_{n\goto +\infty} \int_{y>y_1} \e^2(t_n,y) dy=\alpha_0-{\cal J}.\quad
\ee

Proof of (\ref{masslimit}).\quad 
From the fact that $\e(t_n)\rightharpoonup 0$ in $H^1(\R)$, 
by classical compactness arguments, we have $\e(t_n)\goto 0$ in $L^2_{\rm loc}$ strong.
It follows that
$\overline {\cal J}=\overline{\lim}_{n\goto +\infty} \int_{y<y_1} \e^2(t_n,y) dy$
and $\underline {\cal J}=\underline{\lim}_{n\goto +\infty} \int_{y<y_1} \e^2(t_n,y) dy$
are both independent of $y_1$.

Now, consider $u_n(t,x)=\l^{1/2}(t_n) u(t_n+\l^3(t_n)t,
\l(t_n) x+x(t_n)),$ and for
$t\in [0,T_n)$,
$${\cal J}_{y_1,n}(t)=\int u_n^2(t,x) (1-\psi(x-y_1-{\ts {1\over 4}}x_n(t))  dx.$$
From Lemma \ref{suiteJ}, applied with $t=0$, since $u_n(0)=Q+\e(t_n)$, we have:
$\hbox{for all } \delta_1>0$, $\hbox{for all } y_1<-a_1$, $\exists n_1=n_1(\delta_1,y_1)$, 
such that  $\hbox{for all } n\ge n_1$,
\be\label{valeur}
\left|\int (Q+\e)^2(t_n,y) (1-\psi(y-y_1))dy-{\cal J}\right|\le \delta_1 +C e^{y_1 \over 8}.
\ee
Thus,
$\int_{y<2y_1} \e^2(t_n,y) dy \le {\cal J} + \delta_1+C' e^{y_1 \over 8}.$
Therefore, $\hbox{for all } \delta_1>0$,  for all $y_1<-a_1$,
$${\overline {\cal J}}=\overline{\lim}_{n\goto +\infty} \int_{y<2y_1} \e^2(t_n,y) dy \le {\cal J}+ \delta_1+C' e^{y_1 \over
8}.$$ By $\delta_1\goto 0$, we obtain
${\overline {\cal J}} \le {\cal J} +C' e^{y_1 \over 8},$ and then
by $y_1\goto -\infty$, we have ${\overline {\cal J}} \le {\cal J}.$
Similarly, using (\ref{valeur}), we have ${\underline {\cal J}}\ge {\cal J}.$

By mass conservation applied at $t=t_n$,  for all $y_1\in \R$,
$$\int_{y>y_1} \e^2(t_n)=\alpha_0 - 2 \int \e(t_n) Q -\int_{y<y_1} \e^2(t_n).$$
Therefore, by the fact that $\e(t_n)\rightharpoonup 0$ in $L^2$ weak,
$\lim_{n\goto +\infty} \int_{y>y_1} \e^2(t_n)=\alpha_0-{\cal J},$
and claim (\ref{masslimit}) is proved.

\demo{Step 2.  {\it Mass limit of $\e(t^1_{n'})$ as $n'\goto +\infty$}}
From  the property of almost-monotonicity of the $L^2$ mass on $u_{n'}$,
 we claim that there exist  $a_{n'} \goto +\infty$, and $\delta_{n'}\goto 0$ such that
\be\label{step2}\int_{|y|<a_{n'}} (Q+\e(t^1_{n'}))^2 \le \int Q^2 +\delta_{n'}.\ee
\enddemo

\vglue-20pt
{\it Proof of}\/ (\ref{step2}).   Define
$ {\cal I}_{y_2,n'}(t)=\int u_{n'+1}^2(-t,x) (1-\psi(-x+y_2+{\ts{1\over 4}} x_{n'+1}(-t))) dx.
$
From step 1, there exist $y_{1,n'}\goto -\infty$, $y_{2,n'}\goto +\infty$ and $\delta_{1,n'}\goto 0$,
$\delta_{2,n'}\goto 0$,  such that
$${\cal J}_{y_{1,n'},n'}(0)\ge {\cal J}-\delta_{1,n'},\quad{\rm and}\quad
{\cal I}_{y_{2,n'},n'}(0)\ge (\alpha_0-{\cal J})
-\delta_{2,n'}.$$

By the almost-monotonicity of the $L^2$ mass applied on $u_{n'}$ (see Lemma~\ref{monotonicite}), we have,
for $t^2_{n'}={1\over \l^3(t_{n'})} (t^1_{n'}-t_{n'}),$
$${\cal J}_{y_{1,n'},n'}(t^2_{n'})\ge {\cal J}-\delta_{1,n'}- C e^{y_{1,n'}\over 3}.$$
Therefore, by (\ref{99}), 
$$\int_{y\le {y_{1,n'}}} (Q+\e)^2(t^1_{n'},y) dy
\ge \int_{y\le {\l(t_{n'})\over \l(t^1_{n'})} y_{1,n'}} (Q+\e)^2(t^1_{n'},y) dy
\ge {\cal J}-\delta_{1,n'}- C e^{y_{1,n'}\over 3}.$$

Similarly, we use  ${\cal I}_{y_{2,n'},n'}(t)$, the monotonicity property backwards on $u_{n'+1}$
(in fact, using that $u_{n'+1}(-t,-x)$ is a solution of (\ref{kdv})), and the fact that by (\ref{99})
$$\l(t^1_{n'})\le \l(t_{n'})={(1.1) \l(t_{n'+1})}.$$
We have ${\cal I}_{y_{2,n'},n'}(0)\ge (\alpha_0-{\cal J}) -\delta_{2,n'}$, and so by Lemma \ref{monotonicite},
for $t^3_{n'}=\break{1\over \l^3(t_{n'+1})} (t_{n'+1}-t^1_{n'})$,
$${\cal I}_{y_{2,n'},n'}(t^3_{n'})\ge (\alpha_0-{\cal J}) -\delta_{2,n'}- C e^{-{y_{2,n'}\over 3}}.$$
Therefore,
\begin{eqnarray*}
\int_{y\ge {y_{2,n'}\over 1.1}} (Q+\e)^2(t^1_{n'},y) dy
&\ge& \int_{y\ge {\l(t_{n'+1})\over \l(t^1_{n'})} y_{2,n'}} (Q+\e)^2(t^1_{n'},y) dy
\\
&\ge& (\alpha_0-{\cal J})-\delta_{2,n'}- C e^{-{y_{2,n'}\over 3}}.
\end{eqnarray*}

Now, by  conservation of mass, and preceding estimates:
\bee
\int_{{y_{1,n'}\over 2}\le y\le {y_{2,n'}\over 2}} (Q+\e(t^1_{n'}))^2 &=&
\int u_0^2 -\int_{y\le {y_{1,n'}\over 2}}(Q+\e(t^1_{n'}))^2\\
&& -\ \int_{y\ge {y_{2,n'}\over 2}}(Q+\e(t^1_{n'}))^2\\
&\le & \int u_0^2
-{\cal J} -(\alpha_0-{\cal J}) +\delta_{n'}=\int Q^2 +\delta_{n'},
\eee
where $\delta_{n'}\goto 0$ as $n'\goto +\infty$.
Thus the claim is proved.

\demo{Step 3.  {\it Conclusion}}
It follows from mass and energy properties of $\e(t^1_{n'})$.
We have by (\ref{step2}):
  for all   $a>0$,
$$\overline{\lim}_{n'\goto +\infty}
\int_{|y|<a} (Q+\e(t^1_{n'}))^2 \le \int Q^2\quad {\rm
and~so}\quad \int (Q+\widetilde \e)^2 \le \int Q^2. $$
Therefore, since $\e(t^1_{n'})\rightharpoonup \widetilde \e$ in $H^1$ weak, and by energy arguments as in \S 3.1,
we have the following properties of $\widetilde \e$:
$$
|\widetilde \e|_{H^1}\le C \sqrt{\alpha_0},\quad
\int  yQ_y \widetilde \e=\int  y\left({\ts {Q\over 2}}+yQ_y\right) \widetilde \e=0,
\quad E(Q+\widetilde \e)\le 0.$$

Thus by the variational characterization of $Q$ (see (\ref{charac})), there exist $\l_0>0$ and $x_0\in \R$ such that
$Q+\widetilde \e=\l_0^{1/2} Q(\l_0 x +x_0)$. Then
$|1-\l_0|+|x_0|\le C \sqrt{\alpha_0}$ follows from the smallness of $\widetilde \e$.
From the orthogonality conditions on $\widetilde \e$, we have $\l_0=1$ and $x_0=0$, so that
$\widetilde \e=0$, which is a contradiction.
This concludes the proof of Proposition \ref{wholeseq}. \hfill\qed\enddemo

Therefore Theorem 1 is proved.

\bigbreak\centerline{\bf Appendix A}
\vglue12pt
The objective is to prove the following result
(the function $\Psi_A$ is defined in \S 2.3).

\proclaimtitle{Local Virial relation}
\specialnumber{9}\proclaim{Proposition}\label{APPENDIXA}
There exist $A_0>2${\rm ,} $\alpha_5>0$ and  $\delta_0>0$  such that for $\alpha_0<\alpha_5$
$$\hbox{if} \quad \hbox{for all } s\ge 0,\quad \int yQ_y \e(s)=\int y\left({Q\over 2}+yQ_y\right) \e(s) =0,\quad
\hbox{then}$$
\be\label{ViriallocalappendixA}
 \left(\int \Psi_{A_0} \e^2\right)_s
 \le -\delta_0 \int (\e^2+\e_y^2)e^{-{|y|\over A_0}}
+{1\over \delta_0} \left(\int \e Q \right)^2.
\ee
\endproclaim

\demo{{P}roof} 
Recall the equation satisfied by $\e(s)$:
\bee
 \e_s &=& (L\e)_y   + {\l_s \over \l} \left({Q\over 2} 
      +yQ_y\right) + \left( {x_s \over \l}-1 \right) Q_y  +{\l_s \over \l} \left({\e\over 2} 
      +y\e_y\right) \\
       && +\  \left( {x_s \over \l}-1 \right) \e_y 
         -  (10Q^3\e^2+10Q^2\e^3+5Q \e^4+\e^5)_y,
\eee
where 
$
L \e =-\e_{xx}+\e-5 Q^4 \e.
$
As in \cite{MM1}, proof of Lemma 5, regularization arguments allow us to have the following
relation,  multiplying the equation by $\Psi_A
\e$:
\bee
{1\over 2} \left(\int \Psi_A \e^2\right)_s &=&
-\int L\e \, (\Psi_A \e)_y\\
&& +\, {\l_s\over \l} \int \left({Q\over 2}+yQ_y\right)\Psi_A \e 
   +  \left({x_s\over \l}-1\right) \int Q_y\Psi_A \e\\
&& +\,   {\l_s\over \l} \int \left({\e\over 2}+y \e_y\right)\Psi_A \e 
 +  \left({x_s\over \l}-1\right) \int \e_y \Psi_A \e \\
&&   +\ \int(10Q^3 \e^2+10Q^2\e^3+5Q\e^4)(\Psi_A \e)_y - 5 \int \Psi_A \e^5\e_y.
\eee
Recall that $\Psi_A'=\Phi_A$ and
\bee && \bullet \int L\e \, (\Psi_A \e)_y={3\over 2}\int \e^2_y \Phi_A+{1\over 2}\int \e^2\Phi_A\\
&&\phantom{\bullet \int L\e \, (\Psi_A \e)_y=}
-{1\over 2} \int \e^2 \Phi_A'' -{5\over 2}\int Q^4 \Phi_A \e^2 +10 \int  Q_y Q^3\Psi_A \e^2, \\
&& \bullet \int \left({\e\over 2}+y \e_y\right)\Psi_A \e=-{1\over 2}\int y\Phi_A \e^2, \\
&& \bullet \int \e_y \Psi_A \e=-{1\over 2} \int \Phi_A \e^2,\\
&& \bullet 
\int Q^i \e^{5-i} (\Psi_A \e)_y=
\left({5-i\over 6-i}\right) \int Q^i \Phi_A \e^{6-i} \\
&&\phantom{ \bullet 
\int Q^i \e^{5-i} (\Psi_A \e)_y=}-{i\over 6-i}
\int Q_y Q^{i-1} \Psi_A \e^{6-i}\quad  {\rm for}~i=0,1,2,3.
\eee

By analogy with the notation
\bee
H_{\infty}(\e,\e)
&=& {3\over 2} \int \e_y^2 +{1\over 2}\int \e^2-{5\over 2} \int Q^4 \e^2 +10\int y Q_y Q^3 \e^2,
\eee
we set
$$
H_A(\e,\e)={3\over 2} \int \e^2_y \Phi_A+{1\over 2}\int \e^2 \Phi_A
-{5\over 2}\int Q^4 \Phi_A \e^2 +10\int Q_y Q^3 \Psi_A \e^2.$$

We obtain
\bee
{1\over 2} \left(\int \Psi_A \e^2\right)_s &=&
-H_A(\e,\e) +{\l_s\over \l} \int \left({Q\over 2}+yQ_y\right)\Psi_A \e \\
&&+\ \left({x_s\over \l}-1\right) \int Q_y\Psi_A \e  +{1\over 2} \int \e^2 \Phi_A''
  -  {1\over 2} {\l_s\over \l} \int y\Phi_A \e^2\\
&& -\ {1\over 2} \left({x_s\over \l}-1\right)
\int \Phi_A \e^2  +K_A(\e) +{5\over 6} \int \Phi_A \e^6,
\eee
where 
\bee K_A(\e)&=&{20\over 3} \int Q^3\Phi_A \e^3 -{10} \int Q_y Q^2 \Psi_A \e^3
\\ &&  +\ {15\over 2} \int Q^2 \Phi_A \e^4 -{5} \int Q_y Q \Psi_A  \e^4
+4 \int Q \Phi_A \e^5-\int Q_y \Psi_A \e^5.\eee

The rest of the proof proceeds in four steps. In the first step, we show that under some perturbation of the
orthogonality conditions, $H_{\infty}(\e,\e)$ is still positive. In step 2, we use this result to 
show that under the orthogonality conditions $\int Q \e=\int y\left({Q\over 2}+yQ_y\right) \e=0$, 
for large $A$ (independent of $\alpha_0$), $H_A(\e,\e)$ is positive definite.
In step 3, we remove the condition $\int Q\e =0$ and obtain an additional term in the inequality.
In step 4, we treat all remaining terms in the expression of $\left(\int \Psi_A \e^2\right)_s$
to finish the proof of the proposition.
\enddemo

{\it Step} 1.  Perturbation of the orthogonality conditions in the limit case.
Recall from \cite{MM2} Part B, that there exists $\delta_1>0$ such that 
\be\label{fromliou}
\hbox{if}\quad \int Q \e=\int y\left({Q\over 2}+yQ_y\right) \e=0 \quad
\hbox{then} \quad H_{\infty}(\e,\e)\ge \delta_1 \int \left(\e_y^2 +\e^2\right).
\ee

We claim that
there exists $\delta_2>0$ such that
\be\label{almost}
\hbox{if } |(\e,Q)|+\left| \left(\e, y\left({\ts {Q\over 2}}+yQ_y\right)\right)\right|\le \delta_2 |\e|_{H^1},
\hbox{ then } H_{\infty}(\e,\e)\ge {\delta_1\over 4} \int (\e_y^2+\e^2).\ee

Indeed, take $\e$ satisfying the assumption of (\ref{almost}). Then
$$\e=\e_1+a Q+b y\left({\ts {Q\over 2}}+yQ_y\right)=\e_1+\e_2,$$
where $(\e_1,Q)=\left(\e_1,y\left({\ts {Q\over 2}}+yQ_y\right)\right)=(\e_1,\e_2)=0,$
\be\label{39bis}
a={(\e,Q)\over (Q,Q)},\quad b={\left(\e,y\left({\ts {Q\over 2}}+yQ_y\right)\right)
\over \left(y\left({\ts {Q\over 2}}+yQ_y\right),y\left({\ts {Q\over 2}}+yQ_y\right)\right)}.
\ee
Note that for $\delta_2$ small enough, 
\be\label{comparaison}
{1\over 2} \int (\e_y^2+\e^2) \le \int \e_{1y}^2+\e_1^2 \le 2 \int (\e_y^2+\e^2).
\ee

By bilinearity,
\begin{eqnarray*}
H_{\infty}(\e,\e)&=&H_{\infty}(\e_1,\e_1)+H_\infty(\e_2,\e_2)
+3\int \e_{1y}\e_{2y}\\
&&+\ \int \e_1 \e_2 -5 \int Q^4 \e_1 \e_2 +20 \int y Q_y Q^3 \e_1 \e_2.
\end{eqnarray*}
We have from (\ref{fromliou}), (\ref{comparaison}), (\ref{39bis}),
\bee
H_\infty(\e_1,\e_1)&\ge& \delta_1 \int (\e_{1y}^2+\e_1^2)\ge {\delta_1\over 2}\int (\e_{y}^2+\e^2) , \\ 
H_\infty(\e_2,\e_2)&=&a^2 H_{\infty}(Q,Q)
+b^2 H\left(y\left({\ts {Q\over 2}}+yQ_y\right),y\left({\ts {Q\over 2}}+yQ_y\right)\right)\\
&\le &
C \delta_2^2 \int(\e_y^2+\e^2).
\eee

Now, we have $\int \e_1\e_2=0$ and by integration by parts and (\ref{comparaison}),
$$\left| 3\int \e_{1y}\e_{2y}-5 \int Q^4 \e_1 \e_2 +20 \int y Q_y Q^3 \e_1 \e_2 \right| \le C (|a|+|b|) |\e_1|_{H^1}\le C
\delta_2 |\e|^2_{H^1}.$$

Thus for $\delta_2$ small, independent of $\e$, we have 
$\ds H_\infty(\e,\e)\ge {\delta_1\over 4}  \int (\e_y^2+\e^2).$
This proves the claim.

\demo{Step 2.  {\it Positivity of $H_A(\e,\e)$ for $A$ large}}
We claim: there exists $A_0>2$ such that for $A>A_0$, 
 for all $\e\in H^1$ such that $(\e,Q)=\left(\e, y\left({\ts {Q\over 2}}+yQ_y\right)\right)=0$, 
$$H_A(\e,\e)\ge {\delta_1 \over 8} \int \left(\e_y^2 e^{-{|y|\over A}} +\e^2 e^{-{|y|\over A}}\right).$$

The idea is to write $H_A(\e,\e)$ as $H_{\infty}(\e \sqrt{\Phi_A},\e \sqrt{\Phi_A})$
plus some error terms to be controlled. Then, to conclude, $\e \sqrt{\Phi_A}$ almost satisfies
the orthogonality conditions, in the sense of step 1.

Now,
$$\int \e_y^2 \Phi_A=\int (\e \sqrt{\Phi_A})_y^2 -{1\over 4} \int \e^2 {{\Phi_A'}^2\over \Phi_A}
-\int \e \e_y \Phi'_A,$$
and
$$\int Q_yQ^3 \Psi_A \e^2=
\int y Q_y Q^3 (\e \sqrt{\Phi_A})^2+\int Q_y Q^3 (\Psi_A -y \Phi_A) \e^2.$$
Therefore
\begin{eqnarray*}
H_A(\e,\e)&=&H_{\infty}(\e \sqrt{\Phi_A},\e \sqrt{\Phi_A})\\ &&
-\ {3\over 8} \int \e^2 {{\Phi_A'}^2\over \Phi_A}-{3\over 2} \int \e \e_y \Phi'_A +10 \int Q_y Q^3 (\Psi_A -y \Phi_A) \e^2.
\eee
Now we verify that $\e \sqrt{\Phi_A}$ satisfies (\ref{almost}) for $A$ large:
\begin{eqnarray*}
\left|\int Q \e \sqrt{\Phi_A}\right|&=&\left|\int Q\e +\int Q(\sqrt{\Phi_A}-1)\e\right|\\
&=&\left|\int Q(\sqrt{\Phi_A}-1) \e\right|\le C e^{-{A\over 2}} |\e \sqrt{\Phi_A}|_{L^2},\end{eqnarray*}
since $\Phi_A=1$ on $[-A,A]$, $0<\Phi_A\le 1$ on $\R$, and  the decay property of $Q$.
Similarly,
$|\int y\left({\ts {Q\over 2}}+yQ_y\right)  \e \sqrt{\Phi_A}|=|\int y\left({\ts {Q\over 2}}+yQ_y\right) 
(\sqrt{\Phi_A}-1) \e|\le  C e^{-{A\over 4}} |\e\sqrt{\Phi_A}|_{L^2}.$
Therefore, for $A$ large enough (depending only on $\delta_1$), 
$\e \sqrt{\Phi_A}$ satisfies the condition in (\ref{almost}). Thus
\bee H_{\infty}(\e\sqrt{\Phi_A},\e \sqrt{\Phi_A})&\ge& {\delta_1\over 4} \int (\e \sqrt{\Phi_A})_y^2+\e^2\Phi_A
\\ &=&{\delta_1\over 4} \int (\e_y^2 \Phi_A+\e^2 \Phi_A) 
+{\delta_1\over 16} \int \e^2 {{\Phi'_A}^2\over \Phi_A} +{\delta_1\over 4} \int \e \e_y \Phi'_A,\eee
and so 
\bee H_A(\e ,\e )
&\ge & {\delta_1 \over 4} \int (\e_y^2\Phi_A +\e^2 \Phi_A)
-{3\over 8} \int \e^2 {{\Phi_A'}^2\over \Phi_A}+\left({\delta_1\over 4}-{3\over 2}\right) \int \e \e_y \Phi'_A
\\ && -\ 10 \left|\int Q_y Q^3 (\Psi_A -y \Phi_A) \e^2\right|.
\eee
Note that $0\le |\Phi'_A|\le {C\over A} \Phi_A$, and next, since
$\Psi_A-y \Phi_A=0$ on $[-A,A]$ and $|\Psi_A-y\Phi_A|\le CA,$ we have
$$\left|\int Q_y Q^3 (\Psi_A -y \Phi_A) \e^2\right|\le
C {\rm sup}(e^{-|x|} |\Psi_A- y\Phi_A|) \int Q^3 \e^2
\le C A e^{-A} \int \e^2 e^{-|y|}.$$
Therefore, for $A$ large enough (depending on $\delta_1$), 
$$H_A(\e,\e)\ge {\delta_1\over 8} \int (\e_y^2\Phi_A+\e^2 \Phi_A)\ge {\delta_1\over 8} 
\int \left(\e_y^2 e^{-{|y|\over A}}+\e^2 e^{-{|y|\over A}}\right).
$$
\enddemo

{\it Step} 3. We remove the orthogonality condition $(\e,Q)=0$ using a similar argument to the one of step 1.
Let $\e$ be such that $\left(\e, y\left({\ts {Q\over 2}}+yQ_y\right)\right)=0$.
Let $\e=\e_1+a Q$, where $a=(\e,Q)/(Q,Q)$.
Then
\bee
H_A(\e,\e)&\hskip-7pt =\hskip-7pt& H_A(\e_1,\e_1) +a^2 H_A(Q,Q)+ 3 a \int Q_y \Phi_A \e_{1y} 
+a \int Q\Phi_A \e_1 \\ &\hskip-7pt\hskip-7pt&  -\ 5 a \int Q^5 \Phi_A \e_1
+20 a \int Q_y Q^4 \Psi_A \e_1 \\
&\hskip-7pt\ge\hskip-7pt & {\delta_1 \over 8} \int (\e_{1y}^2 e^{-{|y|\over A}}+\e_1^2 e^{-{|y|\over A}})
\!-\!C a^2 \!-\!C |a| \int (|\e_1|+|\e_{1y}|) e^{-|y|} e^{-{|y|\over 2A}}  \\
&\hskip-7pt\ge\hskip-7pt & {\delta_1 \over 8} \int (\e_{1y}^2 e^{-{|y|\over A}}+\e_1^2 e^{-{|y|\over A}})
\!-\!C a^2 \!-\!C |a| \left(\int (\e_{1y}^2e^{-{|y|\over A}}+\e_1^2  e^{-{|y|\over A}})\right)^{1/2} \\
&\hskip-7pt\ge\hskip-7pt & {\delta_1 \over 16} \int (\e_{1y}^2 e^{-{|y|\over A}}+\e_1^2 e^{-{|y|\over A}})
\!-\!C a^2 \\
&\hskip-7pt\ge\hskip-7pt & {\delta_1 \over 32} \int (\e_y^2 e^{-{|y|\over A}}+\e^2 e^{-{|y|\over A}})
\!-\!C a^2\!-\!C |a| \left(\int (\e_{y}^2e^{-{|y|\over 2A}}+\e^2  e^{-{|y|\over 2 A}})\right)^{1/2} \\
&\hskip-7pt\ge\hskip-7pt & {\delta_1 \over 64} \int (\e_y^2 e^{-{|y|\over A}}+\e^2 e^{-{|y|\over A}})
\!-\!C \left( \e, Q\right)^2.\\
\eee

\demo{Step 4.  {\it Control of the remainding terms}}
Now, we are reduced to showing that by possibly taking a larger $A$ and a smaller $\alpha_0$,
we can control the perturbation  terms in the expression of $\left (\int \Psi_A \e^2\right)_s$
by an arbitrarily small constant multiplied by  $\int (\e_y^2 e^{-{|y|\over A}}+\e^2 e^{-{|y|\over A}})$.
\vglue4pt
First,
$|\Phi_A''(y)|={1\over A^2} |\Phi''({x\over A})|\le {C\over A^2} \Phi_A(y),$ so that
$\left|  \int \e^2\Phi''_A \right| \le {C\over A^2} \int \e^2 \Phi_A.$

Second, we  have by Lemma {\ref{PROPERTIESDECOMPOSITION}}
\bee\left|{\l_s\over \l}\right|+\left|{x_s\over \l}-1\right|&\le& C \left(\int \e^2 e^{-{|y|\over 2}}\right)^{1/2},\\
\noalign{\noindent  and} 
\left|\int \left({Q\over 2}+yQ_y\right) \Psi_A \e\right| &=&\left|\int \left({Q\over
2}+yQ_y\right) (\Psi_A-y) \e\right| \\ &\le&
{\rm sup}(|\Psi_A-y| e^{-{|y|\over 4}}) \int e^{-{|y|\over 4}} |\e|\\& \le&
C e^{-{A\over 8}} \left(\int \e^2 e^{-{|y|\over
2}}\right)^{1/2}.
\eee
Similarly
$\left| \int Q_y \Psi_A \e\right|
\le C e^{-{A\over 8}} \left(\int \e^2 e^{-{|y|\over 2}}\right)^{1/2}.$

Next, since from Lemma \ref{PROPERTIESDECOMPOSITION},  $|\e|_{L^\infty}\le 
|\e|_{H^1}^{1/2} |\e|_{L^2}^{1/2} \le C \sqrt{\alpha_0},$
we have
$$|K_A(\e)|\le C \sqrt{\alpha_0} \int \e^2 \Phi_A,\quad {\rm
and}\quad 
\int \Phi_A \e^6 \le \alpha_0^2 \int \e^2 \Phi_A.$$
\pagebreak

Finally, since $|y| e^{-{|y|\over A}}\le C A e^{-{|y|\over 2A}},$
\bee \left| {\l_s\over \l} \int y\Phi_A \e^2\right|
&\le& C \left(\int e^{-{|y|\over 2}} \e^2\right)^{1/2}
A \int e^{-{|y|\over 2A}} \e^2 \\
&\le& C A \left(\int \e^2\right)^{1/2} \int \e^2 e^{-{|y|\over A}}
\le C\sqrt{\alpha_0} A \int \e^2 e^{-{|y|\over A}}.
\eee

Therefore, we can fix $A_0>2$ large (independent of $\alpha_0$) so that
\begin{eqnarray}
\noalign{\vskip-5pt}\label{pourplustard}&&\\ 
&&\hskip-18pt {1\over 2} \left(\int \Psi_{A_0} \e^2\right)_s
\le -{\delta_1 \over 128} \int (\e_y^2e^{-{|y|\over A_0}}+ \e^2 e^{-{|y|\over A_0}})
+C'' (\e,Q)^2 +\sqrt{\alpha_0} A_0 \int \e^2 e^{-{|y|\over A_0}}.\nonumber
\end{eqnarray}
The constant $A_0>2$ being fixed, we take $\alpha_0>0$ small enough so that
$${1\over 2} \left(\int \Psi_{A_0} \e^2\right)_s
\le -{\delta_1 \over 256} \int (\e_y^2e^{-{|y|\over A_0}} +\e^2 e^{-{|y|\over A_0}})
+C'' (\e,Q)^2.$$ 
Thus the proposition is proved. \hfill\qed\enddemo

\bigbreak\centerline{\bf Appendix B}
\vglue12pt

\demo{{P}roof of Lemma {\rm \ref{suiteJ}}}  
Note first that for all $x_1<0, n\ge 0, t\in [0,T_n)$, $0\le {\cal J}_{x_0,n}(t)\le 2 \int Q^2.$

\demo{Step 1}  We show that for $n\ge n_0$ fixed, the function ${\cal J}_{x_0,n}(t)$ has a limit
${\cal J}_{x_0,n}$ as $t\uparrow T_n$. 
Since $\l_n(t)\le 1,$ by the proof of  Lemma \ref{monotonicite} applied to $u_n$, 
$$\hbox{for all } x_1< -a_1,\hbox{ for all } t\in [0,T_n),\quad {\cal J}'_{x_1,n}(t)\ge -  C x_{nt} e^{-{|x_1|\over 3}} 
e^{-{1\over 4} x_n(t)}.$$
Now,
 for all $x_1< -a_1$, for all $t',t\in [0,T_n),t'<t$,
$$
{\cal J}_{x_1,n}(t) -  {\cal J}_{x_1,n}(t')\ge C
e^{-{|x_1|\over 3}} \left(e^{-{1\over 4}x_n(t)}-e^{-{1\over 4}x_n(t')}\right).$$

Therefore the function $t\mapsto {\cal J}_{x_1,n}(t)-C e^{-{|x_1|\over 3}} e^{-{1\over 4}x_n(t)}$ is nondecreasing
and bounded on $[0,T_n)$. Thus 
it has a limit as $t\uparrow T_n$. Since $e^{-{1\over 4}x_n(t)}$ has a limit as $t\uparrow T_n$
($x_{nt}\ge 0$), this concludes the proof.
\enddemo

\demo{Step 2}
We claim that there exists $C>0$ such that for all $x_1,x_2<0$,  $n_2>n_0$,  $\delta_1>0$,
there exists $n_1=n_1(x_1,x_2,n_2,\delta_1)$ such that
\be\label{claimmono}\hbox{for all } t\in [0,T_{n}),\hbox{ for all } n\ge n_1,\quad 
{\cal J}_{x_1,n}(t)\ge {\cal J}_{x_2,n_2} -\delta_1 -C e^{x_2\over 8}.\hskip.25in \ee
\enddemo

{\it Proof of claim} (\ref{claimmono}). We use the recurrence of $u_n(t)$.
 Let $x_1$, $x_2$, $n_2$ and $\delta_1$  be as in the claim.
Let $n\ge n_1\ge n_2$, with $n_1$    to be chosen later.
We have by direct calculations
$$u_n(t,x)=\l^{1/2}(t_n) u(t_n+\l^3(t_n) t, \l(t_n) x+x(t_n)),$$
and
$$u(t,x)={1\over \l^{1/2}(t_{n_2})} u_{n_2}\left({t-t_{n_2} \over \l^3(t_{n_2})},{x-x(t_{n_2})\over \l(t_{n_2})}\right),$$
so that
$$u_{n}(t,x)={\l^{1/2}(t_n)\over \l^{1/2}(t_{n_2})}
 u_{n_2} \left({\l^3(t_n) \over \l^3(t_{n_2})} t+ {t_n -t_{n_2} \over \l^3(t_{n_2})} ,
{\l(t_n)\over \l(t_{n_2})} x +{x(t_n)-x(t_{n_2})\over \l(t_{n_2})}\right).$$
Note that
$$\l_n( t)={\l(t_n+\l^3(t_n) t)\over \l(t_n)},\quad x_n( t)
={x(t_n+\l^3(t_n) t)-x(t_n)\over \l(t_n)}.$$
By the change of variable $y={\l(t_n)\over \l(t_{n_2})} x +{x(t_n)-x(t_{n_2})\over \l(t_{n_2})}$,
we obtain
\bee
{\cal J}_{x_1,n}(t) & = & \int u_n^2(t,x) \left(1-\psi(x-x_n(t_0)-x_1-{\ts{1\over 4}}(x_n(t)-x_n(t_0)))\right) dx \\
  & = & \int u_{n_2}^2(\hat t_n,y) \left(1-\psi\left({\l(t_{n_2})\over \l(t_n)}(y-y_{1n})\right)\right)dy 
\eee
where 
$\ds \hat t_n={\l^3(t_n)\over \l^3(t_{n_2})} t + {t_n -t_{n_2} \over \l^3(t_{n_2})},$
and
$\ds y_{1n}={x(t_n)-x(t_{n_2}) \over \l(t_{n_2})} +{\l(t_n) \over \l(t_{n_2})} x_n(t_0)
    +{\l(t_n)\over \l(t_{n_2})} x_1 +{1\over 4} {\l(t_n)\over \l(t_{n_2})} (x_n(t)-x_n(t_0)).$
\smallbreak
In view of the preceding
calculations, it is natural to compare ${\cal J}_{x_1,n}(t)$ to ${\cal J}_{x_2,n_2}(\hat t_n)$.
Note that
$${\cal J}_{x_2,n_2}(\hat t_n)=\int u_{n_2}^2(\hat t_n,y) \left(1-\psi(y-y_{2n})\right) dy,$$
where
$\ds y_{2n}=x_{n_2} (t_0) +x_2 +{1\over 4} (x_{n_2}(\hat t_n) -x_{n_2}(t_0)).$
Now,
\bee y_{1n} &=& -{x(t_{n_2}) \over \l(t_{n_2})} +{3\over 4} {x(t_n+\l^3(t_n) t_0) \over \l(t_{n_2})}
+{\l(t_n)\over \l(t_{n_2})} x_1 +{1\over 4} {x(t_n+\l^3(t_n) t) \over \l(t_{n_2})},
\\
y_{2n}&=& {3\over 4} {x(t_{n_2}+\l^3(t_{n_2}) t_0) \over \l(t_{n_2})}  -{x(t_{n_2})\over \l(t_{n_2})}
+x_2 +{1\over 4} {x(t_{n}+\l^3(t_{n}) t) \over \l(t_{n_2})},
\\
y_{1n}-y_{2n}&=&{3\over 4} {x(t_n+\l^3(t_n) t_0)-x(t_{n_2}+\l^3(t_{n_2}) t_0)\over \l(t_{n_2})}
+{\l(t_n)\over \l(t_{n_2})} x_1 -x_2.
\eee
Note that $t_n+\l^3(t_n) t_0\goto T$ as $n\goto +\infty$. Therefore, there exists $n_0=n_0(n_2)$ 
such that 
$x(t_n+\l^3(t_n) t_0)-x(t_{n_2}+\l^3(t_{n_2}) t_0)\ge 0.$
For $n_1$ large depending on $x_1$, $x_2$, $n_2$, we also have ${\l(t_n)\over \l(t_{n_2})} x_1\ge x_2/2.$
Thus
$$y_{1n}-y_{2n}\ge -{x_2\over 2}.$$

Since $\psi'\ge 0$, $x_2<0$, ${\l(t_{n_2})\over \l(t_n)}>1$,  
and  $1-\psi(y-y_{2n})\le C e^{x_2\over 4 \sqrt{3}}\le C e^{x_2\over 8},$ for $y>y_{2n}-{x_2\over 2}$,
we obtain
\bee {\cal J}_{x_1,n}(t) &\ge& \int u_{n_2}^2(\hat t_n,y)
   \left(1-\psi\left({\l(t_{n_2})\over \l(t_n)}(y-y_{2n}+{x_2\over 2})\right)\right)dy \\
&\ge & \int_{y<y_{2n}-{x_2\over 2}} u_{n_2}^2(\hat t_n,y) (1-\psi(y-y_{2n})) dy\\
&\ge & \int u_{n_2}^2(\hat t_n,y) (1-\psi(y-y_{2n})) dy -C e^{x_2\over 8}.
\eee

Therefore,
$${\cal J}_{x_1,n}(t)\ge {\cal J}_{x_2,n_2}(\hat t_n) - C e^{x_2\over 8},$$
where $C$ is a constant.
Now, since
$\hat t_n\goto {T-t_{n_2} \over \l^3(t_{n_2})}=T_{n_2}$, 
it is sufficient to consider a possibly larger  $n_1$  (depending on
$\delta_1$) such that
${\cal J}_{x_2,n_2}(\hat t_n)\ge {\cal J}_{x_2,n_2} -\delta_1.$
Thus claim (\ref{claimmono}) is proved.

\demo{Step 3. {\it Conclusion of the proof}}
We use the claim, with $x_1, x_2<0$.
Passing to the limit $t\goto T_{n}$,  we see that 
$${\cal J}_{x_1,n}\ge {\cal J}_{x_2,n_2} -\delta_1- C e^{x_2\over 8}.$$

Then $\underline{\lim}_{n\goto +\infty} {\cal J}_{x_1,n}\ge {\cal J}_{x_2,n_2} -\delta_1 -C e^{x_1\over 8},$
and next,
$\underline{\lim}_{n\goto +\infty} {\cal J}_{x_1,n}\ge \overline{\lim}_{n\goto +\infty}{\cal J}_{x_2,n}-\delta_1
-C e^{x_1\over 8}.$
Finally, by $\delta_1\goto 0$, we obtain
\be\label{starstar}
\underline{\lim}_{n\goto +\infty} {\cal J}_{x_1,n}\ge \overline{\lim}_{n\goto +\infty}{\cal J}_{x_2,n}
-C e^{x_1\over 8}.
\ee

Note that since for all $t$, ${\cal J}_{x_1,n}(t)$ is 
nondecreasing in $x_1$, the limit ${\cal J}_{x_1,n}$ is also 
nondecreasing in $x_1$, and $\overline{\lim}_{n\goto +\infty} {\cal J}_{x_1,n}$
is still nondecreasing in $x_1$.
We consider 
$${\cal J}=\lim_{x_1\goto -\infty} \left(\overline{\lim}_{n\goto +\infty} {\cal J}_{x_1,n}\right).$$

Note that by (\ref{starstar}), we also have
${\cal J}=\lim_{x_1\goto -\infty} \left(\underline{\lim}_{n\goto +\infty} {\cal J}_{x_1,n}\right).$

First, we prove the lower estimate in (\ref{trois}),
using step 2 and the definition of ${\cal J}$.
By
$\underline{\lim}_{n\goto +\infty} {\cal J}_{x_1,n} \ge {\cal J},$
and (\ref{claimmono}), for $n_2$ large enough depending on $\delta_1$,
we have for all $n\ge n_1(n_2,x_1)$, for all $t\in [0,T_n)$,
\bee {\cal J}_{x_1,n}(t)& \ge&  {\cal J}_{x_1,n_2}-\delta_1 -C  e^{x_1\over 8}\\
&\ge& \underline{\lim}_{n\goto +\infty} {\cal J}_{x_1,n} -2 \delta_1  -C  e^{x_1\over 8} 
\ge {\cal J} -2 \delta_1 - C e^{x_1\over 8}.
\eee

Finally, we prove the upper estimate in (\ref{trois}), by using Lemma \ref{monotonicite}.
Indeed, by this lemma, we have
for all $t, t'\in [0,T_n),t<t',\quad
{\cal J}_{x_1,n}(t)\le {\cal J}_{x_1,n}(t') +C e^{x_1\over 3}.$

Since $t'\goto T_n$, we obtain
for all $t\in [0,T_n),\quad
{\cal J}_{x_1,n}(t)\le {\cal J}_{x_1,n}+C e^{x_1\over 3}.$
Therefore, for $n$ large depending on $\delta_1$:
for all $t\in [0,T_n)$, 
${\cal J}_{x_1,n}(t)\le\break \overline{\lim}_{n\goto +\infty}{\cal J}_{x_1,n}+\delta_1+C e^{x_1\over 3}.$

Note that by (\ref{starstar}), we have 
$\overline{\lim}_{n\goto +\infty}{\cal J}_{x_1,n}\le \underline{\lim}_{n\goto +\infty}{\cal J}_{x_2,n}
+ Ce^{x_1\over 8},$ and so, by $x_2\goto +\infty$:
$\overline{\lim}_{n\goto +\infty}{\cal J}_{x_1,n} \le {\cal J} + C e^{x_1\over 8}.$
Therefore,
$$\hbox{for all } t\in [0,T_n),\quad
{\cal J}_{x_1,n}(t)\le {\cal J} + C e^{x_1\over 8}+\delta_1+C e^{x_1\over 3}.$$
This concludes the proof of the lemma. \hfill\qed
\enddemo\enddemo

\end{document}